\documentclass[aop,MSNbibl,seceqn,citesort,dvips]{arximspdf}
\usepackage{float}

\doi{10.1214/11-AOP714}
\volume{41}
\issue{2}
\pubyear{2013}
\firstpage{871}
\lastpage{913}

\makeatletter
\newtheorem{theorem}{Theorem}[section]
\newtheorem{proposition}{Proposition}[section]
\newtheorem{corollary}{Corollary}[section]
\newtheorem{lemma}{Lemma}[section]
\newproclaim{definition}{Definition}[section]
\newtheorem{notation}{Notational convention}[section]
\newtheorem{claim}{Claim}
\newproclaim{remark}{Remark}[section]

\newcommand\ignore[1]{}

\def\R{\mathbb{R}} 
\def\Z{\mathbb{Z}} 
\def\N{\mathbb{N}} 

\def\sC{\mathcal{C}}
\def\sF{\mathcal{F}}
\def\sL{\mathcal{L}}
\def\sP{\mathcal{P}}

\def\ba{\mathbf{a}}
\def\bb{\mathbf{b}}
\def\bc{\mathbf{c}}

\def\bs{\mathbf{s}}
\def\bu{\mathbf{u}}
\def\bv{\mathbf{v}}
\def\bw{\mathbf{w}}
\def\bx{\mathbf{x}}
\def\by{\mathbf{y}}
\def\bz{\mathbf{z}}

\def\eps{\varepsilon}

\def\ink{\mathsf{ink}}
\def\Fill{\mathsf{Fill}}
\def\hatink{\widehat{\mathsf{ink}}}

\def\IP{\mathsf{IP}}

\def\RW{\mathsf{RW}}
\def\EX{\mathsf{EX}}
\def\fT{\mathsf{T}}
\def\dtv{d_{\mathrm{TV}}}
\def\Unif{\mathrm{Unif}}
\def\bO{\mathbf{O}}
\makeatother

\begin{document}
\begin{frontmatter}

\title{Mixing of the symmetric exclusion processes in terms of the
corresponding single-particle random~walk}
\runtitle{Mixing of symmetric exclusion}

\begin{aug}
\author[A]{\fnms{Roberto Imbuzeiro} \snm{Oliveira}\corref{}\ead[label=e1]{rimfo@impa.br}\thanksref{T1}}
\runauthor{R. I. Oliveira}
\affiliation{IMPA}
\thankstext{T1}{Supported by a \textit{Bolsa de Produtividade em
Pesquisa} and a \textit{Pronex} project on Probability from CNPq, Brazil.}
\address[A]{IMPA\\
Estrada Dona Castorina, 110\\
Rio de Janeiro, RJ 22460-320\\
Brazil\\
\printead{e1}} 
\end{aug}

\received{\smonth{9} \syear{2010}}
\revised{\smonth{8} \syear{2011}}

%
\begin{abstract}
We prove an upper bound for the $\eps$-mixing time of the symmetric
exclusion process on any graph $G$,
with any feasible number of particles. Our estimate is proportional to
${\fT}_{\RW(G)} \ln(|V|/\eps),$
where $|V|$ is the number of vertices in $G$, and ${\fT}_{\RW(G)}$ is the
$1/4$-mixing time of the corresponding single-particle random walk. This
bound implies new results for symmetric exclusion on expanders,
percolation clusters, the giant component of the Erd\"{o}s--R\'{e}nyi
random graph and Poisson point processes in
$\R^d$. Our technical tools include a variant of Morris's chameleon process.
\end{abstract}

%
\begin{keyword}[class=AMS]
\kwd[Primary ]{60J27}
\kwd{60K35}
\kwd[; secondary ]{82C22}.
\end{keyword}

\begin{keyword}
\kwd{Symmetric exclusion}
\kwd{interchange process}
\kwd{mixing time}.
\end{keyword}\vspace*{-3pt}

\end{frontmatter}

\section{Introduction}\label{sec1}
The \textit{symmetric exclusion process} is a continuous-time Markov
chain defined on a \textit{weighted graph} $G=(V,E,\{w_e\}_{e\in E})$,
where $V$ is a set of vertices, $E$ is a set of edges and to each
$e\in E$, we assign a positive weight $w_e>0$. For $k\leq|V|$,
$k$-particle symmetric exclusion on $G$ has the following informal
description.

\textit{Informal description of $\EX(k,G)$}: Start with $k$
indistinguishable particles placed on distinct vertices of $V$. Each
particle moves independently according to the symmetric transition
rates given by the edge weights, except that moves to occupied sites
are suppressed.

This is one of the most basic and best studied processes in the
literature on interacting particle systems
\cite{LiggettIPSBook,LiggettSISBook}. Literally hundreds of papers
have been written on this process, but most of these results apply only
to restricted classes of infinite graphs, such as the lattices $\Z^d$.

Exclusion processes over finite graphs have also been a testbed for the
quantitative
analysis of finite Markov chains. Coupling
\cite{AldousFillRWBook}, comparison arguments
\cite{DiaconisSCComparison}, the martingale method for log-Sobolev
inequalities~\cite{LeeYauLogSobolev,YauLogSobolevExclusion} and
variants of the evolving sets technology\vadjust{\goodbreak}
\cite{MorrisSimpleExclusion,MorrisPeresEvolvingSets} have been
variously applied to this process. Sharp results are known for some
special cases,
such as the complete graph~\cite{LeeYauLogSobolev} and discrete
tori $(\Z/L\Z)^d$
\cite{MorrisSimpleExclusion,YauLogSobolevExclusion}.

In this paper we consider $\EX(k,G)$ over an arbitrary finite graph
and bound its mixing time in terms of the corresponding
single-particle random walk, which we denote by $\RW(G)$. Our result
is very general, but we will see that it nearly matches previously
known mixing results for $\EX(k,G)$ for very specific $G$ and also
gives new results in many interesting classes of examples. We will also
argue that the kind of result presented here
is of conceptual interest.

\subsection{The main result, and why it is interesting}

Recall that the $\eps$-mixing-time of an irreducible continuous-time
Markov chain $Q$ on a finite set $S$,
with transition probabilities $\{q_t(s,s')\}_{s,s'\in S,t\geq0}$,
and stationary (equilibrium) distribution $\pi$, is given by the formula
%
%
\begin{equation}\label{eqdefmixingtime}\fT_Q(\eps)\equiv\inf
\Bigl\{
t\geq0 :  \max_{s\in
S}\dtv(q_t(s,\cdot),\pi)\leq\eps\Bigr\},
\end{equation}
where
$\dtv$ is the total-variation distance; cf. (\ref{eqdtvmeasure}). The
$1/4$ mixing time $\fT_Q(1/4)$ will
also be called \textit{the} mixing time of $Q$. Our main result follows.

\begin{theorem}[(Main result; proven in Section~\protect\ref{seckeysteps})]\label
{thmmain}There exists a universal constant $C>0$ such
for all $\eps\in(0,1/2)$, all connected weighted graphs
$G=(V,E,\{w_e\}_{e\in E})$ with $|V|\geq2$ and all
$k\in\{1,\ldots,|V|-1\}$,
\[
\fT_{\EX(k,G)}(\eps) \leq C  \ln(|V|/\eps)   {\fT}_{\RW(G)}(1/4).
\]
\end{theorem}

Our
bound follows quite naturally if one \textit{assumes} (heuristically)
that the mixing time of $\EX(k,G)$ is not much larger than that of
$k$ independent random walks on~$G$, a process we denote by
$\RW(k,G)$ in what follows:
\[
[\mbox{\textit{Heuristic assumption}}]\quad   \fT_{\EX(k,G)}(\eps)\leq C_0
{\fT
}_{\RW(k,G)}(\eps), \qquad C_0>0\mbox{ universal}.
\]
This assumption, if at all true, is well beyond the reach of present
techniques. However, it
is at least plausible, given that $\RW(k,G)$ and $\EX(k,G)$ are
similar.

It can be shown that ${\fT}_{\RW(k,G)}(\eps)$ and
${\fT}_{\RW(G)}(\eps/k)$ are of the same order if $\eps/k\ll1$; thus
our assumption is
equivalent to
\[
[\mbox{\textit{Heurisitic assumption}}]\quad   {\fT}_{\EX(k,G)}(\eps)\leq
C_1
{\fT}_{\RW(G)}(\eps/k),\qquad C_1>0\mbox{ universal}.
\]
Recall the general inequality ``${\fT}_{\RW(G)}(\delta)\leq C_2
\ln(1/\delta) {\fT}_{\RW(G)}(1/4)$,'' with $C_2>0$ universal, which is
valid for any $0<\delta<1/2$~\cite{AldousFillRWBook}. Applying this
to our assumption, we obtain
\begin{eqnarray}
[\mbox{\textit{Heuristic conclusion}}]\quad {\fT}_{\EX(k,G)}(\eps)\leq
C_3
\ln(k/\eps) {\fT}_{\RW(G)}(1/4),
\nonumber
\\
 \eqntext{C_3>0\mbox{ universal}.}
\end{eqnarray}
Theorem~\ref{thmmain} coincides with this for $k>|V|^c$,
$c>0$ a universal constant; whereas for other $k$ it is a strictly
weaker result.\vadjust{\goodbreak}

We emphasize that what we just presented is \textit{not} a rigorous
proof of Theorem~\ref{thmmain}, since we offer no good grounds for our
heuristic assumption.
What is interesting is that the theorem does give an {a
posteriori} justification for a weakened form of the assumption. We
note that the bound ``${\fT}_{\EX(k,G)}(\eps) \leq C
{\fT}_{\RW(G)}(1/4)\ln(k/\eps)$'' is tight up to constant factors for
some $G$ (e.g., discrete tori $(\Z/L\Z)^d$, $d$ fixed \cite
{MorrisSimpleExclusion}); therefore, in
some sense Theorem~\ref{thmmain} is quite close to the best that one might
hope for.

Many other complex Markov chains are built from simpler
processes that interact; examples appear in, for example, \cite
{AldousFillRWBook,CooperEtAlMultipleRWIPS,MorrisZeroRange}. Given our
main result, it seems reasonable that, at least in some cases, the
mixing time of these complex processes may be bounded in terms of their
constituent parts. Some of the techniques we use to
prove Theorem~\ref{thmmain} are very specific to $\EX(k,G)$, but it
may be that
some of the same ideas will turn out to be useful in other
cases.\looseness=-1

\subsection{Connections with Aldous's conjecture}

Another motivation for our paper is a conjecture of Aldous's for the
\textit{interchange process}, which was recently proved in
\cite{CaputoEtAlInterchange}. The \textit{interchange
process on $G$ with $k\leq|V|$ particles} can be informally
described as follows:

\textit{Informal description of $\IP(k,G)$}: Start with $k$
distinct vertices of $V$ labeled $1,2,\ldots,k$ all remaining vertices
(if any) are labelled ``empty.'' For each edge $e$, switch the labels of
the endpoints of $e$ at rate $w_e$.

One can obtain $\EX(k,G)$ from $\IP(k,G)$ by ``forgetting'' the labels
of the $k$ particles. In particular, the contraction principle \cite
{AldousFillRWBook} implies that ${\fT}_{\EX(k,G)}(\eps)\leq{\fT
}_{\IP
(k,G)}(\eps)$ for all $1\leq k\leq|V|-1$ and all $\eps\in(0,1)$.

Aldous conjectured---and Caputo et al. recently proved
\cite{CaputoEtAlInterchange} (see also~\cite{DiekerInterchange})---that
$\IP(k,G)$ and $\EX(k,G)$ always have the same \textit{spectral
gap} as $\RW(G)$ [or $\RW(k,G)$]. This is a remarkable result, but
it does not say much about the \textit{mixing times} of these
processes, since the bounds for ${\fT}_{\IP(k,G)}(\eps)$ or
${\fT}_{\EX(k,G)}(\eps)$ that can be obtained from the spectral gap are
typically very loose.

Theorem~\ref{thmmain} gives tighter relations
between these mixing times. In the proof of the theorem, we will show
that the bound claimed for $\fT_{\EX(k,G)}(\eps)$ in the theorem
statement in fact holds for $\fT_{\IP(k,G)}(\eps)$ whenever $k\leq
|V|/2$. Our proofs can be adapted to show that
\begin{eqnarray*}
&&\forall\alpha\in(0,1),  \exists C_\alpha>0,  \forall G,  \forall
k\leq\alpha|V| {:}\\
&&\qquad {\fT}_{\IP(k,G)}(\eps)\leq C_\alpha  {\fT}_{\RW
(G)}(1/4)  \ln(|V|/\eps).
\end{eqnarray*}
That is, one can get a bound similar to Theorem~\ref{thmmain} also
for the
interchange process, as long as the fraction of empty sites is bounded
away from $0$. Unfortunately, this leaves out the most interesting case of
$\IP(|V|,G)$, which is a random walk by random transpositions in the
group of permutations of $V$. Fortunately, the restriction on $k$ does
not make a
difference for the exclusion process.

\subsection{Applications and comparison with previous results}

It is not hard to apply Theorem~\ref{thmmain} to specific examples:
all one
needs is a bound for the mixing time of simple random walk on the given
graph, and $\RW(G)$ is typically much easier to analyse than $\EX
(k,G)$. The only example where we know Theorem~\ref{thmmain} gives a
suboptimal
bound is in the case $G=(\Z/L\Z)^d$ with the usual bonds, where the
optimal bound, obtained by Morris~\cite{MorrisSimpleExclusion}, is of
the order $L^2 \ln k$ whereas ours is about $L^2\ln L$ (both for $d$
fixed). Notice that this difference is only relevant for quite small $k$.

%
\begin{table}
\caption{Bounds for $\fT_{\EX(k,G)}(1/4)$ via
Theorem~\protect\ref{thmmain} in examples where no previous bound was
available. We
take $d$ as a fixed parameter and assume $k\approx|V|/2$}\label{table1}
\begin{tabular*}{\textwidth}{@{\extracolsep{\fill}}lc@{}}
\hline
{\textbf{Example}} & {\textbf{Bound for} $\mathbf{\fT}_{\mathbf{\EX}\bolds{(k,G)}}\bolds{(1/4)}$}\\
\hline
$(\Z/L\Z)^d$ with nearest-neighbor bonds~\cite{MorrisSimpleExclusion}&
$|V|^{2/d}\ln|V|$ \\
Typical largest
percolation cluster in $(\Z/L\Z)^d$ \cite
{BenjaminiMosselRWPercolation,PeteConnectivityMixingPercolation} &
$|V|^{2/d}\ln|V|$ \\
Typical Poisson process, in $[0,L]^d$ \cite
{CaputoFaggionatoRWPointProcess} (case $\alpha>d$)& $|V|^{2/d}\ln|V|$
\\
Bounded-degree expanders & $\ln^2|V|$ \\
Giant component of $G_{n,c/n}$, $c>1$~\cite{FountoulakisReedMixingGnp}
& $\ln^3|V|$ \\
\hline
\end{tabular*}
\end{table}

Table~\ref{table1} presents the bounds given by Theorem~\ref
{thmmain} in
examples where no previous bound appears explictly in the literature.
The references are to papers where the mixing times of the
corresponding graphs are computed. We consider only $k\approx|V|/2$
and omit constant factors.

In fairness, we note that a combination of canonical paths, log Sobolev
constants and comparison arguments could in principle be applied to
examples. This method is discussed in Section~\ref{secfaircomparison}
in the
\hyperref[app]{Appendix}. However, we note that:
\begin{itemize}
\item To the best of our knowledge, no good canonical paths bounds have
been worked out for the examples in Table~\ref{table1}, and it might
be hard or impossible to do so;
\item Even if such bounds were obtained, there are natural lower bounds
for how good they can be (cf. Section~\ref{secfaircomparison}), and
Theorem~\ref{thmmain} is at least as good as these lower bounds, up
to the constants
(it is actually better by a $\ln|V|$ factor in the case of expanders).
\end{itemize}

\subsection{Key steps of the proof}\label{seckeysteps}
Our proof of Theorem~\ref{thmmain} can be broken into two main steps.
We first
show that $\IP(2,G)$ always has a mixing time comparable to
$\RW(G)$.\setcounter{footnote}{1}\footnote{Since the single-particle marginal distributions
of $\IP(2,G)$ are given by $\RW(G)$, ${\fT}_{\RW(G)}(1/4)\leq
\fT_{\IP(2,G)}(1/4)$ is immediate from the contraction principle.}

\begin{lemma}\label{lemmixfor2}For any weighted graph $G$,
\[
\fT_{\IP(2,G)}(1/4)\leq20\mbox{,}000  {\fT}_{\RW(G)}(1/4).
\]
\end{lemma}

We then bootstrap the first lemma to a larger number of particles.
\begin{lemma}[(Proven in Section~\protect\ref{secproofkto2})]\label{lemkto2}There
exists a universal constant $K>0$ such that for all connected weighted
graphs $G=(V,E,\{w_e\}_{e\in E})$, all $\eps\in(0,1/2)$ and all
$k\in\{1,\ldots,|V|/2\}$,
\[
{\fT}_{\IP(k,G)}(\eps) \leq K  \fT_{\IP(2,G)}(1/4)  \ln
(|V|/\eps).
\]
\end{lemma}

Before we continue, we show how Theorem~\ref{thmmain} easily follows
from the
two lemmas.

\begin{pf*}{Proof of Theorem~\protect\ref{thmmain}} Combining Lemma~\ref
{lemkto2} with
Lemma~\ref{lemmixfor2} gives
\[
{\fT}_{\IP(k,G)}(\eps)\leq C  {\fT}_{\RW(G)}(1/4) \ln(|V|/\eps
)\qquad\mbox{if }\eps\in(0,1/2) \mbox{ and } k\leq|V|/2
\]
where $C = 20\mbox{,}000 K$. The contraction principle \cite
{AldousFillRWBook} implies
\begin{eqnarray}
{\fT}_{\EX(k,G)}(\eps)\leq{\fT}_{\IP(k,G)}(\eps)\leq C  {\fT
}_{\RW
(G)}(1/4) \ln(|V|/\eps)\nonumber\\
\eqntext{\mbox{if }\eps\in(0,1/2) \mbox{ and }k\leq|V|/2.}
\end{eqnarray}
However, $\EX(k,G)$ and $\EX(|V|-k,G)$ are the same
process with the roles of empty and occupied sites reversed. In
particular, ${\fT}_{\EX(k,G)}(\eps) = \fT_{\EX(|V|-k,G)}(\eps)$
for all
$\eps$.
\end{pf*}

We now give an overview of the main ideas involved in proving the two
lemmas. The proof of Lemma~\ref{lemmixfor2} relies on realizing that there
are two classes of graphs. Some $G$ are ``easy,'' in that two
independent random walkers are likely to meet
by time $O({\fT}_{\RW(G)}(1/4))$ from
\textit{any pair} of initial states. In this case, an argument of Aldous
and Fill's
\cite{AldousFillRWBook} suffices to prove Lemma~\ref{lemmixfor2} (see
Proposition~\ref{propeasy}).

On the other hand, if $G$ is not easy, then for \textit{most} initial
states two independent random walkers are very unlikely to meet by time
$\Omega({\fT}_{\RW(G)}(1/4))$; cf. Proposition~\ref
{propnocollisionfrommost}.
Intuitively, $\IP(2,G)$ and $\RW(2,G)$ are similar in the abscence of
collisions, and we will use this to prove Lemma~\ref{lemmixfor2} over
noneasy graphs. The negative correlation property will be crucial for
this part of the argument; see Remark~\ref{remNEGCORRusedonce} for details.

The proof of Lemma~\ref{lemkto2} is considerably more involved. We
first note
that there are two methods in the literature for moving from mixing of
pairs of particles to many more particles, both of which were
introduced by Morris~\cite{MorrisImproved,MorrisSimpleExclusion}.
The first one~\cite{MorrisImproved} gives bounds for walks on the
symmetric group by random transpositions. Unfortunately, the method
seems to require too much from the process to be useful in our general
setting. Moreover, the bounds given by that method would have a factor
of $\ln(|V|)\ln(1/\eps)$ where ours has a $\ln(|V|/\eps)$ term.\vadjust{\goodbreak}

Morris's other method was introduced in his study of symmetric
exclusion over
$(\Z\setminus L\Z)^d$~\cite{MorrisSimpleExclusion}. The
so-called \textit{chameleon process} features particles that change
color in a way that encodes the conditional distribution of the
$k$th particle in $\IP(k,G)$ given the other $k-1$ particles. It is
this method that we will successfully adapt to prove Lemma~\ref{lemkto2}.

One way to understand Morris's construction is that it reduces the
analysis of mixing to the study of pairwise collisions between
particles. The analysis for $(\Z\setminus L\Z)^d$ is greatly
facilitated by the explicit structure of the graph, something that we
lack in general. This will require certain technical modifications of
Morris's construction, of which we will try to make sense with remarks
in our proofs.

\subsection{Organization}

Section~\ref{secpreliminaries} reviews some preliminary material.
Section~\ref{secdefinitions} discusses $\RW(G)$, $\EX(k,G)$ and
$\IP(k,G)$, presents
their joint graphical construction and reviews the negative correlation
property. Section~\ref{secmixof2} presents the proof of Lemma~\ref
{lemmixfor2}. The
chameleon process is introduced in Section~\ref{secchameleon}. It is
then used
to prove Lemma~\ref{lemkto2} in Section~\ref{secmainproofs}, but
several lemmas are
postponed to Sections~\mbox{\ref{secmiscellany}--\ref
{secdepinkingsarefast}.} It would be pointless to describe these steps
now, but Section~\ref{secremainder} provides an outline of those sections.
Finally, Section~\ref{secconclusion} presents some final remarks, and the
\hyperref[app]{Appendix} contains some technical steps that are not particularly
illuminating.

\section{Preliminaries}\label{secpreliminaries}

\subsection{Basic notation}
$\N=\{0,1,2,3,\ldots\}$ is the set of nonnegative integers and
$\N_+\equiv\N\setminus\{0\}$. For $n\in\N_+$,
$[n]\equiv\{i\in\N_+ :  i\leq n\}=\{1,\ldots,n\}$. If $S$ is a
finite set, $|S|$ is the cardinality of $S$. For any $k\in[|S|]$,
\[
\pmatrix{{S}\vspace*{2pt}\cr{k}} = \{A\subset S :  |A|=k\}
\]
is the set of all size-$k$ subsets of $S$, and
\[
(S)_k = \{\bs=(\bs(1),\ldots,\bs(k))\in S^k \dvtx  \forall i,j\in
[k],
``i\neq j\mbox{''}\Rightarrow``\bs(i)\neq\bs(j)\mbox{''}\}
\]
is the set of all $k$-tuples of distinct elements in $S$.
\begin{notation}The elements of $(S)_k$ will always be denoted by
boldface letters such as $\bx$, with $\bx(i)$ denoting the $i$th
coordinate of $\bx$.\looseness=-1
\end{notation}

Notice
that with these symbols,
\[
\left|\pmatrix{{S}\vspace*{2pt}\cr{k}}\right| = \pmatrix{{|S|}\vspace*{2pt}\cr{k}}, \qquad |(S)_k| =
(|S|)_k.
\]

A \textit{graph} is a couple $H=(V,E)$ where $V\neq\varnothing$ is the
set of \textit{vertices}, and $E\subset{{V}\choose {2}}$ is the set of \textit{edges}.
For each $e\in E$, the two elements $a,b\in V$ such that
$e=\{a,b\}$ are called the \textit{endpoints of $e$}.\vadjust{\goodbreak}

A \textit{weighted graph} is a triple $G=(V,E,\{w_e\}_{e\in E})$, where
$(V,E)$ is a graph, and $w_e>0$, the \textit{weight of edge $e$}, is
positive for each $e\in E$. When a graph $G$ is introduced without
explicitly defining the edge weights, we will assume that they are
all equal to $1$. We will assume throughout this paper that all graphs
we consider are connected.

\subsection{Basic probabilistic concepts}\label{secprobconcepts}

${\sL}[X]$ denotes the \textit{law} or \textit{distribution} of the random
variable $X$.

Given two probability distributions $\mu,\nu$ over the same finite
set $S$, the \textit{total variation distance} between them is given by
several equivalent formulas:
%
%
\renewcommand{\theequation}{\arabic{section}.\arabic{subsection}.\arabic{equation}}
\setcounter{equation}{0}
\begin{eqnarray}\label{eqdtvmeasure}
\dtv(\mu,\nu) &\equiv& \max_{A\subset S}\bigl(\mu(A)-\nu(A)\bigr) \\
\label{eqdtvsupfunctions} &=& \sup_{f:S\to[0,1]}\int f d\mu- \int f  d\nu\\
\label{eqdtvsumpositivepart} &=& \sum_{s\in S}\bigl(\mu(s)-\nu(s)\bigr)_+\\
 &=&\frac{1}{2}\sum_{s\in S}|\mu(s)-\nu(s)|.
\end{eqnarray}
Another equivalent definition of $\dtv$ is
\[
\dtv(\mu,\nu)=\inf\mathbb{P}(X\neq Y),
\]
where the infimum is over all pairs $(X,Y)$ of $S$-valued
random variables with ${\sL}[X]=\mu$ and ${\sL}[Y]=\nu$ [such a
pair is
called a \textit{coupling} of $(\mu,\nu)$]. This implies that for any
pair of $S$-valued random variables $X,Y$ defined over the same
probability space,
\[
\dtv({\sL}[X],{\sL}[Y])\leq\mathbb{P}(X\neq Y).
\]

We will need the following simple fact: if (for $i=1,2$)
$\mu_i,\nu_i$ are probability distributions on the finite set $S_i$,
%
%
\begin{equation}\label{eqdtvproduct}\dtv(\mu_1\times\mu_2,\nu
_1\times
\nu_2)\leq\dtv(\mu_1,\nu_1) + \dtv(\mu_2,\nu_2).
\end{equation}

We will write $\Unif(S)$ for the \textit{uniform distribution} on a set
$S\neq\varnothing$. This is the normalized counting measure on $S$, if
$S$ is finite, or normalized Lebesgue measure over $S$, if $S\subset\R^d$.

\subsection{Markov chains and mixing times}\label{secmixing}

For our purposes it is convenient to define a continous-time Markov
chain over a finite set $S$ as a family of processes
\[
\{(X^s_t)_{t\geq0} : s\in S\}
\]
defined on the same probability space, with the following properties:
\begin{enumerate}[(1)]
\item[(1)] For each $s\in S$, $X^s_0=s$ almost surely.\vadjust{\goodbreak}
\item[(2)]Each $X^s_t$ is a ``c\`{a}dl\`{a}g'' path over $S$: there
exists a divergent sequence
\[
\tau_0=0<\tau_1<\tau_2<\cdots
\]
and a sequence $\{s_i\}_{i\geq0}\subset S$ with $s_0=s$ with
$X^s_t\equiv s_i$ over each interval $[\tau_i,\tau_{i+1})$.
\item[(3)] For each $h\geq0$ and each c\`{a}dl\`{a}g path
$(x_u)_{u\geq0}$ taking values in $S$ [in the sense of (2)],
\[
\mathbb{P}(X^s_{t+h}=s'\vert X^{s}_{t'}=x_{t'}, 0\leq t'\leq t) =
\mathbb{P}(X^{x_t}_h=s')\qquad\mbox{almost surely.}
\]
\end{enumerate}

The last property is the so-called \textit{Markov property}. It also
implies that the law of $(X^{s}_{t+h})_{h\geq0}$ equals that of
$(X^{x_t}_h)_{h\geq0}$ under the above conditioning. It is well known
that any such process is uniquely defined by its transition rates,
\[
q(s,s')\equiv\lim_{\eps\searrow0}\frac{\mathbb{P}(X^s_\eps
=s')}{\eps}\qquad
[(s,s')\in S^2,  s\neq s'],
\]
or equivalenty by its \textit{generator},
\[
Q:f\in\R^S\mapsto Qf(\cdot)\equiv\sum_{s'\in S, s'\neq\cdot
}q(\cdot
,s')\bigl(f(s')-f(\cdot)\bigr).
\]
We will usually make no distinction between a Markov chain and its
generator in our notation.

In this paper we will only work with \textit{irreducible chains}, that is,
chains for which for all $A\subset S$ with $A\neq\varnothing$,
$S\setminus A\neq\varnothing$, there exist $a\in A$, $b\in S\setminus
A$ with $q(a,b)>0$. It is well known that such Markov chains have a
unique \textit{stationary} distribution $\pi$, that is, a distribution
such that if $s_*$ is picked according to $\pi$ independently from the
$(X^s_t)_{t\geq0,s\in S}$, then ${\sL}[X^{s_*}_t]=\pi$ for all
$t\geq
0$. Moreover,
\[
\forall s\in S \qquad \dtv({\sL}[X^s_t],\pi)\searrow0\qquad \mbox{as }t\to
+\infty.
\]
(The symbol ``$\searrow$'' denotes monotone convergence.) The \textit{$\eps
$-mixing time of $Q$} is thus defined as in the \hyperref[sec1]{Introduction},
\[
\fT_Q(\eps)\equiv\inf\Bigl\{t\geq0 : \max_{s\in S}\dtv({\sL
}[X^s_t],\pi
)\leq\eps\Bigr\}\qquad  [\eps\in(0,1)].
\]

We will often need two elementary facts about Markov chains and their
mixing times.

\begin{proposition}[(\cite{LevinPeresWilmerMCBook}, equation (4.36),
page 55)]\label{propmixingpowersof2}Let $Q$ be a Markov chain on finite
state space $S$. Then for all $0<\eps<1/2$,\footnote{The result in
\cite
{LevinPeresWilmerMCBook} is for discrete-time chains, but the proof
trivially extends to continuous time.}
\[
\fT_{Q}(\eps)\leq\lceil\log_2(1/\eps)\rceil  \fT_Q(1/4).
\]
\end{proposition}

\begin{proposition}[(\cite{AldousFillRWBook}, Lemma 7 in Chapter
4)]\label{propmixingseparation}Let $Q$ be a Markov chain on finite
state space $S$ with symmetric transition rates. Then $\pi$ is uniform
over $S$ and
moreover, for all $0<\eps<1/2$ and $t\geq2\fT_{Q}(\eps)$,
\[
\mathbb{P}(X^s_t=s')\geq\frac{(1-2\eps)^2}{|S|},
\]
for all $s,s'\in S$, with the same notation introduced above.
\end{proposition}

We also make the following convenient notational convention.

\begin{notation}\label{notcadlag}By definition, for any c\`{a}dlag
path $(x_t)_{t\geq0}$ there exists a divergent sequence
$t_0=0<t_1<t_2<\cdots$ with $x_t$ constant over $[t_i,t_{i+1})$ for each
$i\geq0$. For $t>0$, we define $x_{t_-}$ to be the state of $x_t$
immediately prior to time $t$. That is,
\[
x_{t_-} \equiv\cases{
x_{t_{i-1}}, & \quad $\mbox{if }t=t_i\mbox{ for some $i\geq
1$;}$\vspace*{2pt}\cr
 x_t, & \quad $\mbox{otherwise.}$}
\]
Notice that $x_{t_-}=x_{t-\delta}$ for all $\delta>0$ sufficiently small.
\end{notation}

\section{Random walks, exclusion and interchange
processes}\label{secdefinitions}

In this section we formally define the main Markov chains in this
paper: $\RW(G),\EX(k,G)$ and $\IP(k,G)$. We
also present the standard \textit{graphical construction} for the three processes
at the same time, and then discuss the negative correlation property
for $\EX(k,G)$. The material in this section is quite classical:
Liggett's books
\cite{LiggettIPSBook,LiggettSISBook} are basic references, and the
manuscript by Aldous and Fill~\cite{AldousFillRWBook} contains some
additional facts on $\IP(k,G)$ as well as a presentation, that is,
somewhat closer in style to ours.

\subsection{Definitions}\label{secdefthreeprocesses}

The three processes we are defined in terms of the same \textit{weighted graph} $G=(V,E,\{w_e\}_{e\in E})$ with $V$ finite; cf.
Section~\ref{secpreliminaries}. We will be implicitly assuming that
$G$ is connected,
in which case one can easily show that the chains defined below are
irreducible. It will be useful to define the transpositions
\[
f_e:x\in V\mapsto
\cases{
b, & \quad $\mbox{if }x=a,$\vspace*{2pt}\cr
a, & \quad $\mbox{if }x=b,$\vspace*{2pt}\cr
x, & \quad $\mbox{otherwise}.$}
\]
We also write $f_e(A) = \{f_e(a) :  a\in A\}$ and $f_e(\bx) =
(f_e(\bx
(i)))_{i=1}^k$ for $A\in{{V}\choose{k}}$ and $\bx\in(V)_k$ (resp.).

\textit{Simple random walk on $G$}, denoted by $\RW(G)$, is the
continuous-time Markov chain with state space $V$ and transition rates
\[
q(u,v) \equiv\cases{
w_e,  &\quad$\mbox{if }f_e(u)=v;$\vspace*{2pt}\cr
0, &\quad $\mbox{otherwise}$}\qquad[(u,v)\in(V)_2].
\]
We will also consider the process $\RW(k,G)$ that corresponds to $k$
such random walks performed simultaneously and independently. Since the
transition rates of these process are also symmetric, it follows that
the stationary distribution of $\RW(k,G)$ is $\Unif(V^k)$ for all
$k\in
\N_+$.

\textit{The $k$-particle symmetric exclusion process on $G$}, denoted
by $\EX(k,G)$, is the continuous-time Markov chain with state space
${{V}\choose{k}}$ and
transition rates
\[
q^{\{k\}}(A,B) \equiv\cases{
w_e, & \quad $\mbox{if }f_e(A) = B;$\vspace*{2pt}\cr
0, &\quad $\mbox{otherwise}$}\qquad
\left[(A,B)\in
\left(\pmatrix{{V}\vspace*{2pt}\cr{k}}\right)_2\right].
\]
The transition rates are again symmetric, and the stationary
distribution of $\EX(k,G)$ is $\Unif({{V}\choose{k}})$.

\textit{The $k$-particle interchange process on $G$}, denoted by
$\IP(k,G)$, has state space $(V)_k$. The transition rates of $\IP(k,G)$
are given by
\[
q^{(k)}(\bx,\by) \equiv\cases{
w_e, & \quad $\mbox{if }f_e(\bx) =\by;$\vspace*{2pt}\cr
0, &\quad $\mbox{otherwise}$}\qquad [(\bx,\by)\in
((V)_k)_2].
\]
This process also has symmetric transition rates, and its stationary
distribution is $\Unif((V)_k)$.

\subsection{The standard graphical construction}\label{secstandard}

We now present the standard \textit{graphical construction} of these three
processes. Graphical constructions are standard tools in the study of
interacting particle systems~\cite{LiggettIPSBook} and are usually
attributed to Harris in the literature. The basic construction
presented here will be elaborated upon later in the paper; see
Section~\ref{secchameleon}. For brevity, we omit all proofs in this
subsection.

Set $W=\sum_{e\in E}w_e$. We need a marked Poisson process, that is, a
pair of independent ingredients given as follows:

\begin{enumerate}[(1)]
\item[(1)] A Poisson process $\sP=\{\tau_1\leq\tau_2\leq\tau
_3\leq
\cdots\}\subset[0,+\infty)$ with rate
$W$.
\item[(2)] An i.i.d. sequence of $E$-valued random variables (``markings'')
$\{e_n\}_{n\in\N}$, with
\[
\forall n\in\N\qquad
\mathbb{P}(e_n=e)=w_e/W.
\]
\end{enumerate}

Let $0\leq t\leq s<+\infty$ be given. We define a random permutation
$I_{(t,s]}:V\to V$ associated with the time interval $(t,s]$ as
follows: if $\sP\cap(t,s] =\varnothing$, $I_{(t,s]}$ is the identity map
on $V$. If, on the other hand,
\[
\sP\cap(t,s]\equiv\{\tau_j : m\leq j\leq n\}\neq\varnothing,
\]
we set $I_{(t,s]} = f_{e_n}\circ f_{e_{n-1}}\circ\cdots\circ
f_{e_m}$; that is, $I_{(t,s]}$ is the composition of each transposition
$f_{e_j}$ corresponding to $\tau_j\in(t,s]$, and the transpositions
are composed in the order they appear. We also set $I_t\equiv
I_{(0,t]}$ for $t>0$ and $I_{(t,t]}={}$identity map over~$V$.

\begin{remark} Strictly speaking, we should worry about what happens
if $\sP\cap(t,s]$ is infinite, or (more generally) some finite
interval $(a,b]$ in $[0,+\infty)$ has infinite intersection with $\sP$.
However, since the probability of any of this holding is $0$, we will
simply ignore these issues.
\end{remark}

Notice the following simple properties:

\begin{proposition}[(Proof omitted)]\label{propsemigroup}For all $0\leq
t\leq s\leq r$, $I_{(t,r]} = \break I_{(s,r]}\circ I_{(t,s]}$.
\end{proposition}
\begin{proposition}[(Proof omitted)]\label{propbackwards}For all $0\leq
t\leq s<+\infty$,\break ${\sL}[I_{(t,s]}]={\sL}[I_{(t,s]}^{-1}]$.
\end{proposition}
\begin{proposition}[(Proof omitted)]\label{propIndependent}Let $0\leq
t_0<t_1<t_2<\cdots<t_k$. Then the maps $I_{(t_{i-1},t_i]}$, $1\leq
i\leq
k$, are independent.
\end{proposition}

\begin{notation}\label{notIt}We ``lift'' the
random maps $I_{(t,s]}$ to permutations of ${{V}\choose{k}}$ and $(V)_k$,
which we also denote by $I_{(t,s]}$:
\begin{eqnarray*}
I_{(t,s]}(A)&\equiv&\bigl\{I_{(t,s]}(a) :  a\in A\bigr\} \qquad\left [A\in\pmatrix{{V}\vspace*{2pt}\cr{k}}\right],\\
I_{(t,s]}(\bx)&\equiv&\bigl(I_{(t,s]}(\bx(1)),I_{(t,s]}(\bx(2)),\ldots,I_{(t,s]}(\bx(k))\bigr)
\qquad
[\bx\in(V)_{k}].
\end{eqnarray*}
For brevity, we will often write $x^I_t,A^I_t,\bx^I_t$ instead of
$I_t(x),I_t(A),I_t(\bx)$ (resp.).
\end{notation}

The key property of the graphical construction follows:

\begin{proposition}[(Proof omitted)]\label{propgraphical}Let $t_0\geq
0$. Then:
\begin{enumerate}[(1)]
\item[(1)] For each $x\in V$, the process $\{I_{(t_0,t+t_0]}(x)\}
_{t\geq0}$ is a realization of $\RW(G)$ with initial state $x$.
\item[(2)] For each $A\in{{V}\choose{k}}$, the process $\{
I_{(t_0,t+t_0]}(A)\}_{t\geq0}$ is a realization of $\EX(k,G)$ with
initial state $A$.
\item[(3)] For each $\mathbf{x}\in(V)_k$, the process $\{
I_{(t_0,t+t_0]}(\mathbf{x})\}_{t\geq0}$ is a realization of $\IP(k,G)$
with initial state $\mathbf{x}$.
\end{enumerate}
\end{proposition}

\subsection{The negative correlation property}

$\EX(k,G)$ enjoys important \textit{negative correlation properties}. In
this paper we only need a very special result, which is contained in
any of
\cite
{LiggettIPSBook,LiggettInvariantMeasures2,AndjelCorrelationSymmetricExclusion}.

\begin{lemma}\label{lemNEGCORR1}Given $A\in{{V}\choose{k}}$, let $\{
A^I_t\}_{t\geq0}$ be
a realization of $\EX(k,G)$ starting from $A$. Then for all $\bu\in
(V)_2$---that is, for all distinct $\bu(1),\bu(2)\in V$---we have
\[
\mathbb{P}\bigl(\{\bu(1)\in A^I_t\}\cap\{\bu(2)\in A^I_t\}\bigr)\leq\mathbb
{P}\bigl(\bu(1)\in A^I_t\bigr) \mathbb{P}\bigl(\bu(2)\in A^I_t\bigr).
\]
\end{lemma}

Using the construction in the previous section, we can write the above
inequality as
\[
\mathbb{P}\bigl(\{I_t^{-1}(\bu(1))\in A\}\cap\{I_t^{-1}(\bu(2))\in A\}
\bigr)\leq\mathbb{P}\bigl(I_t^{-1}(\bu(1))\in A\bigr) \mathbb{P}\bigl(I_t^{-1}(\bu
(2))\in A\bigr).
\]
The following is then immediate from Proposition~\ref{propbackwards}.

\begin{corollary}[(Proof omitted)]\label{corNEGCORR}Given $\bu\in(V)_2$,
let $\{\bu^I_t\}_{t\geq0}$ be
a realization of $\IP(2,G)$ starting from $\bu$. Then for all
$A\subset V$,
\[
\mathbb{P}\bigl(\{\bu^I_t(1)\in A\}\cap\{\bu^I_t(2)\in A\}\bigr)\leq\mathbb
{P}\bigl(\bu ^I_t(1)\in A\bigr) \mathbb{P}\bigl(\bu^I_t(2)\in A\bigr).
\]
\end{corollary}

\section{The dynamics of pairs of particles}\label{secmixof2}

The goal of this section is to prove Lemma~\ref{lemmixfor2}. We fix a
weighted graph $G=(V,E,\{w_e\}_{e\in E})$ for the remainder of the
section (and of the paper). The definitions of
$\RW(G)$, $\RW(k,G)$, $\EX(k,G)$ and $\IP(k,G)$ are all relative to this
graph.

\subsection{Some facts on $\RW(2,G)$ and $\IP(2,G)$}

Much of this section will involve comparisons between $\IP(2,G)$ and
$\RW(2,G)$. The following notational convention will be useful.

\begin{notation}\label{notx^R}Given $\bx\in V^2$, $\{\bx
^{R}_t\equiv
(\bx^R_t(1),\bx^R_t(2)) : t\geq0\}$ denotes a realization of $\RW
(2,G)$ from initial state $\bx$. That is, the trajectories of $\bx
_t^R(1),\bx_t^R(2)$ are independent realizations of $\RW(G)$ with
respective initial states $\bx(1),\bx(2)$.
\end{notation}

We collect several simple facts about $\RW(2,G)$ and $\IP(2,G)$ that
we will need later on. The first one is obvious, for example, from the
graphical construction.

\begin{proposition}[(Proof omitted)]\label{propmarginalsRWIP}For
$i=1,2$, ${\sL}[\bx^R_t(i)]={\sL}[\bx^I_t(i)]$.
\end{proposition}

The next proposition is a direct consequence of (\ref{eqdtvproduct}).
\begin{proposition}[(Proof omitted)]The mixing times of $\RW(2,G)$
satisfy~ ${\fT}_{\RW(2,G)}(\eps)\leq{\fT}_{\RW(G)}(\eps/2).$
\end{proposition}

\begin{proposition}\label{propdecreasesinRW2}Let $k\in\N$ be given.
Then ${\fT}_{\RW(2,G)}(2^{-k})\leq(k+1)\times\break {\fT}_{\RW(G)}(1/4).$
\end{proposition}
\begin{pf} Follows from the previous proposition combined with
Proposition~\ref{propmixingpowersof2}.
\end{pf}

The next lemma has the following meaning. Suppose $t$ is so large that
$\bx^R_t$ is close to equilibrium. In this case, $\mathbb{E}
[\phi(\bx^R_t)]$
is close to the uniform average of $\phi$ over $V^2$, for all mappings
$0\leq\phi\leq1$. The lemma shows that $\mathbb{E}[\phi(\bx
^I_t)]$ cannot be
much larger than that average. This will require the negative
correlation property; cf. Corollary~\ref{corNEGCORR}.\vadjust{\goodbreak}
\begin{lemma}\label{lemsymmetric}Let $\phi:V^2\to[0,1]$. Then
\[
\forall\eps\in(0,1/16),  \forall t\geq\fT_{\RW(G)}(\eps),
\forall
\bx\in(V)_2 \qquad \mathbb{E}[\phi(\bx^I_t)]\leq8\sqrt
{\eps} + 9 \sum_{\bv\in
V^2}\frac{\phi(\bv)}{|V|^2}.
\]
\end{lemma}
\begin{pf}Define the ``good set'' of all $a\in V$ with nearly uniform
probability
\[
\mathrm{Good}\equiv\biggl\{a\in V {:}  \max_{i=1,2}\biggl|\mathbb
{P}\bigl(\bx_t^R(i)=a\bigr)
- \frac{1}{|V|}\biggr|\leq\frac{2\sqrt{\eps}}{|V|}\biggr\}.
\]
We will show toward the end of the proof that
%
%
\begin{equation}\label{eqprobbadset}\mathbb{P}(\bx^I_t\notin
\mathrm{Good}^2)\leq8\sqrt{\eps},
\end{equation}
which (since $0\leq\phi\leq1$) implies
%
%
\begin{equation}\label{eqprobbadset2}\mathbb{E}\bigl[\phi(\bx
^I_t){\mathbb{I}_{(V)_2\setminus\mathrm{Good}^2}}(\bx^I_t)\bigr]\leq
8\sqrt
{\eps}.
\end{equation}
On the other hand, notice that
\begin{eqnarray*}
\mathbb{E}[\phi(\bx^I_t){\mathbb{I}_{\mathrm{Good}^2}}(\bx^I_t)] &=& \sum
_{\ba\in(\mathrm{Good}^2)\cap(V)_2}\mathbb{P}\bigl(\bx_t^I=(\ba(1),\ba
(2))\bigr)
\phi(\ba
)\\
&\leq& \sum_{\ba\in(\mathrm{Good}^2)\cap(V)_2}\mathbb{P}\Biggl(\bigcap
_{i=1}^2\bigl\{ \bx _t^I(i)\in\{\ba(1),\ba(2)\}\bigr\}\Biggr) \phi(\ba)\\
\mbox{(Cor.~\ref{corNEGCORR})}&\leq& \sum_{\ba\in(\mathrm{Good}^2)\cap(V)_2}\prod
_{i=1}^2\mathbb{P}\bigl(\bigl\{\bx_t^I(i)\in\{\ba(1),\ba(2)\}\bigr\}\bigr) \phi(\ba
)\\
\mbox{(Prop.~\ref{propmarginalsRWIP})}&=& \sum_{\ba\in(\mathrm{Good}^2)\cap
(V)_2}\prod_{i=1}^2\mathbb{P}\bigl(\bigl\{\bx_t^R(i)\in\{\ba(1),\ba(2)\}\bigr\}
\bigr) \phi
(\ba)\\
\mbox{[$\ba(i)\in\mathrm{Good}$]}&\leq& \sum_{\ba\in(\mathrm{Good}^2)\cap(V)_2}\biggl(\frac{2+4\sqrt{\eps}}{|V|}\biggr)^2
\phi(\ba
)\\
\bigl(\sqrt{\eps}\leq1/4\bigr)&\leq& 9 \sum_{\ba\in
V^2}\frac{\phi(\ba)}{|V|^2}.
\end{eqnarray*}
Combining this with (\ref{eqprobbadset2}) finishes the proof, except for
(\ref{eqprobbadset}). To prove that, we let $\mathrm{Bad} = V\setminus
\mathrm{Good}$. Notice that
\[
\frac{\sqrt{\eps} |\mathrm{Bad}|}{|V|}\leq\sum_{a\in V}\frac
{1}{2}\biggl\{
\biggl|\mathbb{P}\bigl(\bx_t^R(1)=a\bigr)-\frac{1}{|V|}\biggr|+\biggl|\mathbb
{P}\bigl(\bx _t^R(2)=a\bigr)-\frac{1}{|V|}\biggl|\biggr\}
\]
as each $a\in\mathrm{Bad}$ contributes at least $\sqrt{\eps}/|V|$ to the
sum. But the RHS equals
\[
\dtv({\sL}[\bx^R_t(1)],\Unif(V))+\dtv({\sL}[\bx^R_t(2)],\Unif
(V))\leq
2\eps
\]
since $t\geq\fT_{\RW(G)}(\eps)$. We deduce
\[
\frac{\sqrt{\eps} |\mathrm{Bad}|}{|V|}\leq2\eps\quad\mbox{or equivalently
}\quad|\mathrm{Good}|\geq\bigl(1-2\sqrt{\eps}\bigr) |V|.
\]
Moreover, $\mathbb{P}(\bx_t^R(i)=a)\geq(1-2\sqrt{\eps})|V|^{-1}$
for all
$a\in\mathrm{Good}$, hence
\[
\mathbb{P}\bigl(\bx_t^R(i)\in\mathrm{Good}\bigr)\geq\frac{|\mathrm{Good}|}{|V|}
\bigl(1-2\sqrt
{\eps}\bigr)\geq\bigl(1-2\sqrt{\eps}\bigr)^2\geq1 - 4\sqrt{\eps}.
\]
Inequality (\ref{eqprobbadset}) now follows from
\begin{eqnarray*}
\mathbb{P}(\bx^I_t\notin\mathrm{Good}^2)&\leq
&\mathbb{P}\bigl(\bx ^I_t(1)\notin\mathrm{Good}\bigr) + \mathbb{P}\bigl(\bx
^I_t(2)\notin\mathrm{Good}\bigr)\\
\mbox{(Proposition~\ref{propmarginalsRWIP})}&=&\mathbb{P}\bigl(\bx^R_t(1)\notin\mathrm{Good}\bigr) +
\mathbb{P}\bigl(\bx^R_t(2)\notin\mathrm{Good}\bigr)\leq8\sqrt{\eps}.
\end{eqnarray*}
\upqed\end{pf}

\subsection{When collisions are nearly as fast as
mixing}\label{seceasymeeting}

Recalling Notational convention~\ref{notx^R}, we define the first
\textit{meeting time}
$M(\bx)$ of $\RW(2,G)$ started from $\bx\in V^2$ as the smallest
$t_0\geq0$ such that $\bx_{t_0}^R(1)=\bx^R_{t_0}(2)$ (this is a.s.
finite by ergodicity). We will also write
\[
M_{\geq t}(\bx)=\inf\{h_0\geq0 :  \bx_{t+h_0}^R(1)=\bx
^R_{t+h_0}(2)\}
\]
for the time until the first meeting after $t$ (this is a
``time-shifted'' meeting time).

The following definition will be crucial for our analysis.

\begin{definition}\label{defeasy}We say that a weighted graph $G$ is
easy if
\[
\sup_{\bx\in V^2}\mathbb{P}\bigl(M(\bx)> 20\mbox{,}000 {\fT}_{\RW
(G)}(1/4)\bigr)\leq1/8.
\]
\end{definition}

We note that all long enough paths and cycles are examples of easy
graphs. Noneasy graphs include $(\Z/L\Z)^d$ for $d\geq2$ fixed and
$L$ sufficiently large, as well as large expander graphs. The next
proposition proves Lemma~\ref{lemmixfor2} for all easy graphs via a coupling
argument due to Aldous
and Fill.

\begin{proposition}\label{propeasy}Lemma~\ref{lemmixfor2} holds for
all easy
weighted graphs.
\end{proposition}
\begin{pf*}{Proof sketch} Given $G$, Aldous and Fill [\cite*{AldousFillRWBook}, Chapter~14, Section 5]
construct a coupling of $\IP(|V|,G)$
started from two different states $\bu,\bv$. Letting $\{\bu^I_t,\bv
^I_t\}_{t\geq0}$ denote the coupled trajectories, the following
property holds: for each $1\leq i\leq|V|$, $\bu^I_t(i),\bv^I_t(i)$
behave as independent random walks up to their first meeting time,
which we denote by $M_i$. After this time $M_i$,\vadjust{\goodbreak} $\bu^I_t(i)=\bv
^I_t(i)$, that is, the two processes move together. This implies
\begin{eqnarray*}
\forall t\geq0\qquad  \dtv({\sL}[\bu^I_t],{\sL}[\bv^I_t])&\leq&\mathbb
{P}(\bu ^I_t\neq\bv^I_t)\leq\sum_{i=1}^{|V|}\mathbb{P}\bigl(\bu
^I_t(i)\neq\bv ^I_t(i)\bigr)\\
&\leq&
\sum_{i=1}^{|V|}\mathbb{P}(M_i>t).
\end{eqnarray*}
It is easy to adapt this to a coupling of $\IP(2,G)$ starting from
given $\bx,\by\in(V)_2$, so that, if
$\{\bx^I_t,\by^I_t\}_{t\geq0}$ denotes the coupled trajectories, we have
\[
\forall t\geq0 \qquad \dtv({\sL}[\bx^I_t],{\sL}[\by^I_t])\leq\mathbb
{P}(M_1>t) +
\mathbb{P}(M_2>t).
\]
Now both $M_1$ and $M_2$ are the meeting times of
independent random walkers on $G$, which shows that
\[
\forall t\geq0\qquad \sup_{\bx,\by\in(V)_2}\dtv({\sL}[\bx^I_t],{\sL
}[\by ^I_t])\leq2\sup_{\bz\in V^2}\mathbb{P}\bigl(M(\bz)> t\bigr).
\]
For $t={20\mbox{,}000 {\fT}_{\RW(G)}(1/4)}$ and $G$ easy, the RHS is $\leq
1/4$. By convexity, this implies that
\[
\sup_{\bx\in(V)_2}\dtv\bigl({\sL}[\bx^I_t],\mathrm{Unif}((V)_2)\bigr)\leq
\frac{1}{4}.
\]
In other words, $\fT_{\IP(2,G)}(1/4)\leq
20\mbox{,}000 {\fT}_{\RW(G)}(1/4)$.
\end{pf*}

\begin{remark}\label{remaldousfill} Aldous and Fill's argument
actually proves Theorem~\ref{thmmain} for
all easy graphs; see~\cite{AldousFillRWBook}, Chapter 14, Section 5
for details.
\end{remark}

\subsection{Long time to meet in noneasy graphs}

We now consider what happens when $\IP(2,G)$ is performed on a graph,
that is, not easy. Our first goal is to show that independent random
walkers take a relatively long time to meet from \textit{most} initial
states in $V$.

\begin{proposition}\label{propnocollisionfrommost}Assume $G=(V,E,\{
w_e\}_{e\in E})$ is not easy. Then
\[
\frac{1}{|V|^2}\sum_{\bv\in V^2}\mathbb{P}\bigl(M(\bv)\leq20{\fT
}_{\RW (G)}(1/4)\bigr)\leq\frac{1}{125}.
\]
\end{proposition}

\begin{remark} In general we cannot guarantee that $\mathbb
{P}(M(\bv )<20\times\break{\fT }_{\RW(G)} (1/4))$ is uniformly small over all
$\bv\in(V)_2$. In
particular, the probability of collision from adjacent $\bv(1),\bv(2)$
might be much greater than the above bound.
\end{remark}

\begin{pf*}{Proof of the Proposition} Set $T={\fT}_{\RW(G)}(1/4)$.
Since $G$ is not easy, there exists some $\bx\in V^2$ with
%
%
\setcounter{equation}{0}
\begin{equation}\label{eqisbigpacas}\mathbb{P}\bigl(M(\bx)>20\mbox{,}000 T\bigr)>1/8.\vadjust{\goodbreak}
\end{equation}
Consider some $k\in\N$. Using the Markov property and the notation
introduced in Section~\ref{seceasymeeting}, one can write
\[
\mathbb{P}\bigl(M(\bx)>40kT\bigr) =
\mathbb{E}\bigl[{\mathbb{I}_{\{M(\bx)>40(k-1)T\}}}\mathbb
{P}\bigl(M_{\geq 40(k-1)T}(\bx)>40T\vert\bx^R_{40(k-1)T}\bigr)\bigr].
\]
The conditional probability in the RHS equals $\mathbb{P}(M(\by
)>40T)$ for
$\by
=\bx^R_{40(k-1)T}$, hence
\begin{eqnarray*}
\mathbb{P}\bigl(M(\bx)>40kT\bigr)&\leq& \Bigl(\sup_{\by\in
V^2}\mathbb{P}\bigl(M(\by)>40T\bigr)\Bigr) \mathbb{P}\bigl(M(\bx)>40(k-1)T\bigr)\\
\mbox{(\ldots induction\ldots)}&\leq& \Bigl(\sup_{\by\in
V^2}\mathbb{P}\bigl(M(\by)>40T\bigr)\Bigr)^{k}.
\end{eqnarray*}
Applying this to
$k=500$ and using the bound in (\ref{eqisbigpacas}) gives the
following with room to spare:
\[
\sup_{\by\in V^2}\mathbb{P}\bigl(M(\by)>40T\bigr)\geq8^{-1/500}\geq
e^{-3/500}\geq
\frac
{497}{500}.
\]
Fix some $\by\in V^2$ achieving this supremum. Notice that
$M(\by)>40T$ holds if and only if $\by^R_t(1)\neq\by^R_t(2)$ for
all $0\leq t\leq40T$. If, that is, the case, $\by^R_{20T+h}(1)\neq
\by^R_{20T+h}(2)$ for all $0\leq h\leq20T$. Using the Markov
property as before, we see that
\begin{eqnarray*}
\frac{497}{500} &\leq&\mathbb{P}\bigl(M(\by)>40T\bigr)\leq\mathbb{P}\bigl(M_{\geq
20T}(\by )>20T\bigr)\\
&=&\sum
_{\bv\in V^2}\mathbb{P}(\by^R_{20T}=\bv)\mathbb{P}\bigl(M(\bv)>20T\bigr).
\end{eqnarray*}
Moreover, by (\ref{eqdtvsupfunctions}),
\begin{eqnarray*}
&&\sum_{\bv\in V^2}\mathbb{P}(\by^R_{20T}=\bv)\mathbb{P}\bigl(M(\bv
)>20T\bigr)\\
&&\qquad\leq\sum
_{\bv\in
V^2}\frac{\mathbb{P}(M(\bv)>20T)}{|V|^2} + \dtv({\sL}[\by
^R_{20T}],\mathrm{Unif}(V^2)).
\end{eqnarray*}
Hence
\[
\frac{497}{500} - \dtv({\sL}[\by^R_{20T}],\mathrm{Unif}(V^2))\leq\sum
_{\bv
\in V^2}\frac{\mathbb{P}(M(\bv)> 20T)}{|V|^2}.
\]
We finish by noting that, by Proposition~\ref{propdecreasesinRW2},
$20T\geq
{\fT}_{\RW(2,G)}(2^{-19})$, hence
\[
\dtv({\sL}[\by^R_{20T}],\mathrm{Unif}(V^2))\leq2^{-19}\leq\frac{1}{500},
\]
and therefore
\[
\sum_{\bv\in V^2}\frac{\mathbb{P}(M(\bv)>20T)}{|V|^2}\geq\frac
{496}{500} = 1
- \frac{1}{125}.
\]
\upqed\end{pf*}

\subsection{If meeting takes a long time, $\IP(2,G)$ and $\RW(2,G)$
are similar}

We have just shown that the meeting is unlikely to be smaller than
$20{\fT}_{\RW(G)}(1/4)$ from most initial states.
We now show that $\IP(2,G)$ is similar to $\RW(2,G)$ until the first
meeting time.

\begin{proposition}\label{proplocal1}For any $\bx\in(V)_2$ and
$s\geq0$,
\[
\dtv({\sL}[\bx^R_s],{\sL}[\bx^I_s])\leq\mathbb{P}\bigl(M(\bx)\leq s\bigr).
\]
\end{proposition}

We will only need the following simple corollary (proof omitted) in
what follows.
\begin{corollary}\label{corcoupleRW}For any $\bx,\by\in(V)_2$ and
$s\geq0$,
\[
\dtv({\sL}[\bx^I_s],{\sL}[\by^I_s])\leq\mathbb{P}\bigl(M(\bx)\leq s\bigr)
+\mathbb{P}\bigl(M(\by)\leq s\bigr) + \dtv({\sL}[\bx^R_s],{\sL}[\by^R_s]).
\]
\end{corollary}

\begin{pf*}{Proof of Proposition~\protect\ref{proplocal1}} We present a
coupling of $\{
\bx^I_t\}_{t\geq0}$ and $\{\bx^R_t\}_{t\geq0}$ such that the two
processes agree up to $M(\bx)$. The proposition then follows from the
coupling characterization of $\dtv(\cdot,\cdot\cdot)$; cf.
Section~\ref{secprobconcepts}.

Our coupling is given by a continuous-times Markov chain on
$S=(V)_2\times V^2$ with transition rates given by $q(\cdot,\cdot
\cdot
)$. The state space can be split into two parts, $\Delta\equiv\{(\bz
,\bz) : \bz\in(V)_2\}$ and its complement $\Delta^c$.

\begin{itemize}
\item\textit{Transition rule} 1: The transition rates from any pair $(\bx
,\by)\in\Delta^c$ to any other pair in $S$ are the same as those of
independent realizations of $\RW(2,G)$ and $\IP(2,G)$.
\item\textit{Transition rule} 2: The transition rates from a pair $(\bx
,\bx)\in\Delta$ are determined as follows:
\begin{itemize}
\item\textit{Transition rule} 2.1: For each $e\in E$ with $|e\cap\{\bx
(1),\bx(2)\}|=1$,
\[
q((\bx,\bx),(f_e(\bx),f_e(\bx)))=w_e;
\]
\item\textit{Transition rule} 2.2: If $e\in E$ satisfies $e=\{\bx(1),\bx
(2)\}$,
\[
\cases{
((\bx,\bx),(f_e(\bx),(\bx(1),\bx(1)))) =
w_e,\vspace*{2pt}\cr
q((\bx,\bx),(\bx,(\bx(2),\bx(2)))) = w_e.}
\]
\item\textit{Transition rule} 2.3: All other potential transitions have
rate $0$.
\end{itemize}
\end{itemize}
Inspection of the marginals reveals that this indeed gives a coupling
of $\{\bx^R_t\}_{t\geq0}$ and $\{\bx^I_t\}_{t\geq0}$ when started
from an initial state $(\bx,\bx)\in\Delta$. Moreover, the two
processes can only differ after a transition has occurred according to
rule $2.2$. The first time when this happens is precisely the first
meeting time of $\{\bx^R_t\}_{t\geq0}$.
\end{pf*}

\subsection{Proof of the mixing time bound for $\IP(2,G)$}\label
{secproofmixfor2}

We now use the tools developed above in order to prove Lemma~\ref{lemmixfor2}.

\begin{pf*}{Proof of Lemma~\protect\ref{lemmixfor2}} The case of easy
graphs is
covered by Proposition~\ref{propeasy}, so assume $G=(V,E,\{w_e\}
_{e\in E})$ is not
easy. Let $\bx,\by$ be given and $T\equiv{\fT}_{\RW(G)}(1/4)$. Notice
that for all $A\subset(V)_2$, if $\{\bx_{t}^I\}_{t\geq0}$, $\{\by
_{t}^I\}_{t\geq0}$ are defined over the same probability space,
\begin{eqnarray*}\mathbb{P}(\bx_{40T}^I\in A) - \mathbb{P}(\by
_{40T}^I\in A)& =&
\mathbb{E}[\mathbb{P}(\bx_{40T}^I\in A\vert\bx
^I_{20T})-\mathbb{P}(\by _{40T}^I\in A\vert\by^I_{20T})]\\
&\leq&\mathbb{E}\bigl[\dtv
\bigl(\mathbb{P}(\bx_{40T}^I\in\cdot\vert\bx^I_{20T}),\mathbb{P}(\by
_{40T}^I\in \cdot\vert\by^I_{20T})\bigr)\bigr].
\end{eqnarray*}
Maximizing over $A$ yields
%
%
\setcounter{equation}{0}
\begin{eqnarray}\label{eqlocalcondpunchline}
&&\dtv({\sL}[\bx
_{40T}^I],{\sL}[\by_{40T}^I])
\nonumber
\\[-8pt]
\\[-8pt]
\nonumber
&&\qquad\leq\mathbb{E}\bigl[\dtv\bigl(\mathbb
{P}(\bx _{40T}^I\in\cdot\vert\bx^I_{20T}),\mathbb{P}(\by
_{40T}^I\in\cdot\vert \by^I_{20T})\bigr)\bigr].
\end{eqnarray}
By the Markov property and Corollary~\ref{corcoupleRW},
\begin{eqnarray*}
&&\dtv\bigl(\mathbb{P}(\bx_{40T}^I\in\cdot\vert\bx
^I_{20T}=\bv),\mathbb{P}(\by_{40T}^I\in\cdot\vert\by^I_{20T}=\bw
)\bigr)\\
&&\qquad =\dtv({\sL}[\bv ^I_{20T}],{\sL}[\bw^I_{20T}])\\
&&\qquad\leq\mathbb
{P}\bigl(M(\bv)\leq20T\bigr) +\mathbb{P}\bigl(M(\bw)\leq20T\bigr) +
\dtv
({\sL}[\bv^R_{20T}],{\sL}[\bw^R_{20T}]).
\end{eqnarray*}
Proposition~\ref{propdecreasesinRW2} implies the third term in the
RHS is $\leq
2^{-19}$ for any $\bv,\bw$. Using
this in conjunction with (\ref{eqlocalcondpunchline}), we obtain
%
%
\begin{equation}\label{eqlocalcondpunchline2}\dtv({\sL}[\bx
_{40T}^I],{\sL}[\by_{40T}^I])\leq\mathbb{E}[\phi(\bx
^I_{20T})] + \mathbb{E}[\phi(\by^I_{20T})] + 2^{-19},
\end{equation}
where $\phi(\bz) = \mathbb{P}(M(\bz)\leq20T)$. Notice that $0\leq
\phi
\leq
1$. We may apply Lemma~\ref{lemsymmetric} and the fact that $20T\geq
{\fT}_{\RW(G)}(2^{-20})$ (cf. Proposition~\ref
{propmixingpowersof2}) to deduce
%
%
\begin{equation}\label{eqcrucialNEGCORR}\mathbb{E}[\phi(\bx
^I_{20T})]\leq
2^{-7} + 9\sum_{\bv\in V^2}\frac{\mathbb{P}(M(\bv)\leq20T)}{|V|^2}.
\end{equation}
Applying the same reasoning to $\phi(\by^I_{20T})$ and plugging the
results into (\ref{eqlocalcondpunchline2}), we obtain
%
%
\begin{equation}\label{eqlocalcondpunchline3}\quad\dtv({\sL}[\bx
_{40T}^I],{\sL}[\by_{40T}^I])\leq18\sum_{\bv\in V^2}\frac{\mathbb
{P}(M(\bv)\leq 20T)}{|V|^2} + 2^{-6} + 2^{-19}.
\end{equation}
Finally, we use the
fact that $G$ is not easy, combined with
Proposition~\ref{propnocollisionfrommost}, to deduce
%
%
\begin{equation}\label{eqlocalcondpunchline33}\dtv({\sL}[\bx
_{40T}^I],{\sL}[\by_{40T}^I])\leq\frac{18}{125}+ 2^{-9} +
2^{-6}\leq1/4
\end{equation}
with room to spare. By convexity,
%
%
\begin{equation}\label{eqlocalcondpunchline4}\dtv({\sL}[\bx
_{40T}^I],\mathrm{Unif}((V)_2))\leq1/4.
\end{equation}
Since this holds for all $\bx\in(V)_2$, we have
$\fT_{\IP(2,G)}(1/4)\leq40T$, which implies Lemma~\ref{lemmixfor2} for
noneasy graphs.
\end{pf*}

\begin{remark}\label{remNEGCORRusedonce} The first inequality in
(\ref{eqcrucialNEGCORR}) follows from Lemma~\ref{lemsymmetric},
which is a
consequence of the negative correlation property; cf. Lemma~\ref{lemNEGCORR1}
and Corollary~\ref{corNEGCORR}. This is the first crucial use we make of
negative correlation in this paper.
\end{remark}

\section{The chameleon process}\label{secchameleon}

In the previous section we determined the order of magnitude of the
mixing time of $\IP(2,G)$. Going beyond two particles will require an
important additional idea, that is, based on
Morris's paper~\cite{MorrisSimpleExclusion}. His idea is to
introduce the so-called \textit{chameleon process} to keep track of the
conditional distribution of one particle in $\IP(k,G)$. We will need
a different process, which will nevertheless call by the same name.

\subsection{A modified graphical construction}\label{secmodgraphical}

We will need consider a variant of the construction of $\IP(k,G)$
presented in Section~\ref{secstandard}. Consider \textit{three} independent
ingredients:
\begin{enumerate}[(1)]
\item[(1)] A Poisson process $\sP=\{\tau_1\leq\tau_2\leq\tau
_3\leq\cdots\}\subset[0,+\infty)$ with rate
$2W$.
\item[(2)] An i.i.d. sequence of $E$-valued random variables
$\{e_n\}_{n\in\N}$, with $\mathbb{P}(e_n=e)=w_e/W$.
\item[(3)] An i.i.d. sequence of \textit{coin
flips} $\{c_n\}_{n\in\N}$ with $\mathbb{P}(c_n=1)=\mathbb{P}(c_n=0)=1/2$.
\end{enumerate}

Recall the definition of $f_e$ from Section~\ref
{secdefthreeprocesses}, and set
$f_e^1=f_e$, $f_e^0={}$the identity
function. We modify the definition of the maps $I_{(t,s]}$ from
Section~\ref{secstandard} as follows: if $\sP\cap(t,s]=\varnothing$,
$I_{(t,s]}$ is the identity map, as before. Otherwise,
\[
\sP\cap(t,s] = \{\tau_{n}<\tau_{n+1}<\cdots<\tau_{m}\},
\]
and we set
\[
I_{(t,s]} = f^{c_m}_{e_m}\circ\cdots\circ
f^{c_{n+1}}_{e_{n+1}}\circ f^{c_{n}}_{e_{n}}.
\]
The \textit{thinning property} of the Poisson process implies that $\{\tau
_n :  c_n=1\}$ is a Poisson process with rate $W$. One can use this
to show that:
\begin{proposition}[(Proof omitted)]\label{propmodisthesame}The joint
distribution of the maps $I_{(t,s]}$, $0\leq t<s<+\infty$, is the same
as in Section~\ref{secstandard}.
\end{proposition}

\subsection{The chameleon process}\label{secdefchameleon}

The \textit{chameleon process} is built on top of the modified graphical
construction. The definition of the process will depend on a parameter
$T>0$ which we call the \textit{phase length}, for reasons that will
become clear later on.

Given ${\by}\in(V)_{k-1}$, let $\bO(\by)\equiv
\{\by(1),\ldots,\by(k-1)\}$ denote the set of vertices that ``occupied''
by the coordinates of ${\by}$. The chameleon process will
be a continuous-time, time-inhomogeneous Markov chain with state
space
%
%
\begin{eqnarray}\label{eqdefstatespacechameleon}\sC_k(V)&\equiv&\{
({\bz
},R,P,W) :  \bz\in(V)_{k-1};
\nonumber
\\[-8pt]
\\[-8pt]
\nonumber
&&\hspace*{4pt} \mbox{the sets $\bO(\bz),R,P,W$
partition }V\}.
\end{eqnarray}
Notice that we do allow any of the $R,P,W$ to be empty in the above
definition. For a given $(\bz,R,P,W)\in\sC_k(V)$, it will be convenient
to refer to the vertices in the sets $\bO(\bz),R,P,W$ as \textit{black},
\textit{red}, \textit{pink} and \textit{white} (resp.). Notice that any vertex $v\in V$
will belong to one of these color classes.

The evolution of the process from initial state $({\bz},R,P,W)$ will be
denoted by
\[
\{({\bz}^C_t,R^C_t,P^C_t,W^C_t)\}_{t\geq0}.
\]
By definition, this process will only be updated at the times $\tau_n$
($n\in\N$) given by the Poisson process and at deterministic times
$2iT$, $i\in\N$. Moreover, the updates at times $\tau_n$ are of
different kinds depending on whether $\tau_n\in((2i-2)T,(2i-1)T]$ for
some $i\in\N_+$, or $\tau_n\in((2i-1)T,2iT]$ for some $i\in\N_+$.
Finally, we will \textit{define} for convenience,
\[
(\bz^C_{0_-},R^C_{0_-},P^C_{0_-},W^C_{0_-})=({\bz},R,P,W)
\]
and will allow an instantaneous change at time $t=0$: that is,
\[
\mbox{it might happen that }(\bz^C_{0},R^C_{0},P^C_{0},W^C_{0})\neq
({\bz},R,P,W).
\]
The three update rules are described in Box~\ref{boxchameleon}.
\begin{remark} Technically, this process is \textit{not} c\`{a}dl\`
{a}g, as it changes at time $0$. We will nevertheless continue to use
$t_-$ (cf. Notational convention~\ref{notcadlag}) with the proviso
for $t=0$ that we have
just described.
\end{remark}

\floatstyle{ruled}
\newfloat{floatbox}{thp}{lob}[section]
\floatname{floatbox}{Box}
\begin{floatbox}[t]\caption{The three kinds of
updates in the chameleon process.}\label{boxchameleon}
\begin{itemize}
\item\textsf{Constant-color phases}: If $t=\tau_n\in((2i-2)T,(2i-1)T]$
for some $i\in\N_+$, update
%
%
\begin{equation}\label{equpdateconstantcolor}({\bz
}^C_{t},R^C_{t},P^C_{t},W^C_{t})=(f^{c_n}_{e_n}({\bz
}^C_{t_-}),f^{c_n}_{e_n}(R^C_{t_-}),f^{c_n}_{e_n}(P^C_{t_-}),f^{c_n}_{e_n}(W^C_{t_-})).
\end{equation}
That is, the states of the
endpoints of $e_n$ are flipped if $c_n=1$, and nothing happens if $c_n=0$.
\item
\textsf{Color-changing phases.} If $t=\tau_n\in((2i-1)T,2i
T]$, for $i\in\N_+$, update as above \textit{unless}:
\begin{enumerate}[1.]
\item[1.]$e_n=\{w,r\}$ has a white endpoint $w\in W^C_{t_-}$ and a red
endpoint $r\in R^C_{t_-}$;
\item[2.]
$|P^C_{t_-}|<\min\{|R^C_{t_-}|,|W^C_{t_-}|\}$.
\end{enumerate}
If $(1)$ and $(2)$ hold, $r$
and $w$ both become pink, and we call $t$ a pinkening time.
%
%
\begin{equation}\label{equpdatecolorchange}({\bz
}^C_{t},R^C_{t},P^C_{t},W^C_{t})= ({\bz}^C_{t_-},R^C_{t_-}\setminus\{
r\},P^C_{t_-}\cup\{r,w\},W^C_{t_-}\setminus\{w\}).
\end{equation}
\item
\textsf{Depinking times.} If $t=2iT$ with $i\in\N$ ($t=0$ or
$t$ lies at the end of a color-changing phase) \textit{and}
$|P^C_{t_-}|\geq\min\{|W^C_{t_-}|,|R^C_{t_-}|\}$ (more pink than
either white or red), flip a fair coin $d_{i}$, and make all pink
particles become red or white depending on whether $d_i$ comes out
heads or tails (resp.).
%
%
\begin{equation}\label{equpdatedepinking}({\bz
}^C_{t},R^C_{t},P^C_{t},W^C_{t})=\cases{
({\bz}^C_{t_-},R^C_{t_-}\cup P^C_{t_-},\varnothing
,W^C_{t_-}),& \quad $d_i=1;$\vspace*{2pt}\cr ({\bz}^C_{t_-},R^C_{t_-},\varnothing
,W^C_{t_-}\cup
P^C_{t_-}), & \quad $d_i=0.$}
\end{equation}
\ignore{Note that $2kT$ might be a depinking time if
$|P^C_{2kT_-}|=\min
\{|W^C_{2kT_-}|,|R^C_{2kT_-}|\}=0$}
\end{itemize}
\end{floatbox}

\begin{remark} We briefly note that our chamaleon process is more
complicated than Morris's process~\cite{MorrisSimpleExclusion}. In
brief: his process does not have constant-color phases and will depink
right when the number of pink particles exceeds the minimum of red and
white. The second difference is a matter of convenience, but the first
one will be fundamental at key steps of our argument.
\end{remark}
%
\subsection{Two basic properties}\label{secbasicpropertieschameleon}

The next two results will be useful later on. We only sketch the proofs.

\begin{lemma}\label{lemdiscreteMarkov}Let
\[
(\hat{\bz}_i,\hat{R}_i,\hat{P}_i,\hat{W}_i) = \mbox{the value of
}(\bz
^C_{2iT_-},R^C_{2iT_-},P^C_{2iT_-},W^C_{2iT_-}) \qquad (i\in\N).
\]
Then $\{(\hat{\bz}_i,\hat{R}_i,\hat{P}_i,\hat{W}_i)\}_{i\in\N}$
is a
discrete-time, time-homogeneous Markov chain. Moreover, if $D_j$ is the
$j$th depinking time of the process, then $\hat{D}_j\equiv D_j/2T$ is a
stopping-time for this discrete-time Markov chain.
\end{lemma}
\begin{pf*}{Proof Sketch} Markovianity and time-homogeneity are
obvious. To prove the stopping time property, it suffices to check that
(setting $D_0=0$),
\[
\forall j>0 \qquad \frac{D_j}{2T} = \inf\biggl\{i>\frac{D_{j-1}}{2T} :
|\hat{P}_{i}|\geq\min\{|\hat{R}_{i}|,|\hat{W}_{i}|\}\biggr\},
\]
where we allow the $\inf$ to be $+\infty$ if the set is empty or
$D_{j-1}=+\infty$.
\end{pf*}

\begin{lemma}\label{lemjustlikeI}Suppose $(\bz
^C_{2iT},R^C_{2iT},P^C_{2iT},W^C_{2iT})$ is the state of the chameleon
process at time $2iT$ (i.e., at the beginning of a constant-color
phase). Then
\[
\bigl(\bz^C_{(2i+1)T},R^C_{(2i+1)T},P^C_{(2i+1)T},W^C_{(2i+1)T}\bigr) = (I(\bz
^C_{2iT}),I(R^C_{2iT}),I(P^C_{2iT}),I(W^C_{2iT})),
\]
where $I=I_{(2iT,(2i+1)T]}$ is the map defined in the modified
graphical construction.
\end{lemma}
\begin{pf}By inspection.
\end{pf}

\subsection{The chameleon process and conditional distributions}\label
{secchameleonisrelated}
We now explain the relationship between the chameleon process and
conditional distributions.

\begin{notation}\label{notlastcoordinate}$\bx=(\bx(1),\ldots,\bx
(k))\in
(V)_k$ is represented as a
pair $(\bz,x)$, where $\bz=(\bx(1),\ldots,\bx(k-1))\in(V)_{k-1}$ and
$x=\bx(k)\in V\setminus\mathbf{O}(\bz)$. [Notice that $\bx^I_t =
(\bz^I_t,x^I_t)$ for all $t\geq0$.]
\end{notation}

\begin{proposition}[(Proof omitted)]\label{propchameleonisrelated}Given
an initial state $\mathbf{x}=(\mathbf{z},x)\in(V)_k$ for $\IP(k,G)$, set
$R=\{x\}$, $P=\varnothing$ and $W=V\setminus(\mathbf{O}(\bz)\cup\{x\})$.
Consider the interchange process $\{\mathbf{x}^I_t=(\mathbf{z}_t^I,x^I_t)\}_{t\geq0}$ started from state $\mathbf{x}$ and the
chameleon process $\{(\mathbf{z}^C_t,R^C_t,P^C_t,W_t^C)\}_{t\geq0}$
started from configuration $(\mathbf{z},R,\break P,W)\in\sC_k(V)$. Then
%
%
\setcounter{equation}{0}
\begin{equation}\label{eqchameleonisrelated}\forall t\geq0,
\forall
\mathbf{b}=(\mathbf{c},b)\in(V)_k\qquad  \mathbb{P}(\mathbf{x}^I_t=\mathbf{b}) =
\mathbb
{E}\bigl[\ink_t(b){\mathbb{I}_{\{\mathbf{z}^C_t=\mathbf{c}\}}}\bigr],
\end{equation}
where
%
%
\begin{equation}\label{eqdefrhoink}\ink_t(v)\equiv{\mathbb{I}_{\{
v\in R^C_t\} }} + \frac{{\mathbb{I}_{\{v\in P^C_t\}}}}{2}\qquad (v\in V).
\end{equation}
\end{proposition}

This is almost identical (up to changes in notation) to \cite[Lemma
1]{MorrisSimpleExclusion}, and we omit its proof. It will be useful to
think of $\ink_t(v)$ as the amount of ``red ink'' at vertex $v\in V$: a
red vertex has one unit of red ink, a pink vertex has half a unit, and
black or white vertices have no ink. We will see below that the \textit{total} amount of red ink in the system determines the rate of
convergence to equilibrium of $\IP(k,G)$.

\section{From $2$ to $k$ particles via the chameleon process}\label
{secmainproofs}

In this section we present the proof of Lemma~\ref{lemkto2}, modulo several
lemmas about the chameleon process that we will prove later. We then
outline the remainder of the paper.

\subsection{\texorpdfstring{Proof of Lemma~\protect\ref{lemkto2}}{Proof of Lemma 1.2}}\label{secproofkto2}
\mbox{}
\begin{pf}We assume we have defined a chameleon process over $\sC
_k(V)$ as in Section~\ref{secdefchameleon}. We will take the notation and
definitions from that section for granted. We also define

%
\begin{equation}\label{eqdeftotalink}\ink_t\equiv\sum_{v\in V}\ink
_t(v) = |R^C_t| + \frac{|P^C_t|}{2}\qquad  (t\geq0).
\end{equation}
We note for later reference that
%
%
\begin{equation}\label{eqtotalink2}\ink_t\equiv\sum_{v\in
V\setminus
\mathbf{O}(\bz_t^I)}\ink_t(v)
\end{equation}
since the vertices in $\mathbf{O}(\bz^I_t)$ have zero red ink.\vadjust{\goodbreak}

We have argued in Proposition~\ref{propchameleonisrelated} that the
distribution of
$\IP(k,G)$ started from $\bx=(\bz,x)\in(V)_k$ corresponds to a
chameleon process started from $(\bz,\{x\},\varnothing,V\setminus(\bO
(\bz
)\cup\{x\}))$. Letting $\ink_t^\bx$ denote the value of $\ink_t$ in
that chameleon process, we will show that:

\begin{lemma}[(Proven in Section~\protect\ref{secproofchameleonmixing})]\label
{lemchameleonmixing}The following inequality holds for all $1\leq
k\leq|V|-1$:
\[
\sup_{\bx\in(V)_k}\dtv\bigl({\sL}[\bx_t^I],\Unif((V)_k)\bigr)\leq2k \sup
_{\bx
\in(V)_k}\mathbb{E}\biggl[1 - \frac{\ink^{\bx
}_t}{|V|-k+1} \Big| \Fill\biggr]
\]
where
\[
\Fill\equiv\Bigl\{\lim_{t\to+\infty}\ink_t^\bx=|V|-k+1\Bigr\}.
\]
\end{lemma}

The main goal is to bound the expected value in the RHS of the
inequality in Lemma~\ref{lemchameleonmixing}. Fix some $\bx\in
(V)_k$, and let
$D_j(\bx)$ denote the $j$th depinking time for the chameleon process
corresponding to $\bx$. Also set $\hatink^\bx_j\equiv\ink^\bx
_{D_j(\bx
)}$ for this process. We will show in Proposition~\ref
{propinfmanydepinkings} that
there are infinitely many depinking times, that is, there are
infinitely many times of the form $2iT$ at which the number of pink
particles is at least as large as the minimum of the numbers of white
and red. The definition of the chameleon process implies that $\ink
_t^\bx$ can only change at depinking times, hence for any $t\geq0$
$\ink^\bx_t=1$ if $t<D_1(\bx)$ and $\ink^\bx_t = \hatink_{j}^\bx
$ if
$D_j(\bx)\leq t<D_{j+1}(\bx)$ for some $j$. We deduce that
\begin{eqnarray*}
1 - \frac{\ink^{\bx}_t}{|V|-k+1}&\leq&\sup_{m\geq
j}\biggl(1 - \frac{\hatink^\bx_m}{|V|-k+1}\biggr) + {\mathbb{I}_{\{
D_j(\bx )>t\} }}\\
&\leq&\sum_{m\geq j}\biggl(1 - \frac{\hatink^\bx
_m}{|V|-k+1}\biggr) +
{\mathbb{I}_{\{D_j(\bx)>t\}}}.
\end{eqnarray*}
Taking expectations, we see that the RHS of the inequality in
Lemma~\ref{lemchameleonmixing} is at most
%
%
\begin{equation}\label{eqgettingthere111}
2k \sup_{\bx\in
(V)_k}\biggl\{
\sum_{m\geq j}\mathbb{E}\biggl[1 - \frac{\hatink^\bx
_m}{|V|-k+1}\Big| \Fill\biggr]
+ \mathbb{P}\bigl(D_j(\bx)\geq t\vert \Fill\bigr)\biggr\}.
\end{equation}
A simple (but technical) proposition proven in the \hyperref[app]{Appendix} will
take care of the first term.
\begin{proposition}[(Proven in Section~\ref
{secproofdiscreteinkdecay})]\label
{propdiscreteinkdecay}For all $\ell\geq1$ and $\bx\in(V)_{k}$,
\[
\mathbb{E}\biggl[1-\frac{\hatink^{\bx}_\ell}{|V|-k+1}\Big\vert\Fill
\biggr]\leq\sqrt{|V|-k+1}
\biggl(\frac{71}{72}\biggr)^\ell.
\]
\end{proposition}
We thus have
\begin{eqnarray}\label{eqgettingthere112}&& 2k \sup_{\bx\in
(V)_k}\biggl\{
\sum_{m\geq j}\mathbb{E}\biggl[1 - \frac{\hatink^\bx
_m}{|V|-k+1}\Big|\Fill\biggr]
+ \mathbb{P}\bigl(D_j(\bx)\geq t\vert\Fill\bigr)\biggr\}
\nonumber
\\[-8pt]
\\[-8pt]
\nonumber
&&\qquad\leq C_2
|V|^{3/2}
e^{-c_1j} + 10k\sup_{\bx\in(V)_k}\mathbb{P}\bigl(D_j(\bx)\geq t\vert
\Fill\bigr),
\end{eqnarray}
where $c_1 = \ln(72/71)>0$ and $C_2=720$ are universal constants.

Bounding $\mathbb{P}(D_j(\bx)\geq t\vert\Fill)$ is the key step in
the proof.
Up to now all of our results have been valid for all values of $k,|V|$
and of the phase length parameter $T>0$. The next lemma will require
restrictions on these values.
\begin{lemma}[(Proven in Section~\ref{secproofexpmomentDj})]\label
{lemexpmomentDj}There exist universal constants $C_3,C_4>0$, such that
if $|V|\geq300$, $T\geq C_3 \fT_{\IP(2,G)}(1/4)$ and $k/|V|\leq
1/2$, then
\[
\forall\bx\in(V)_k,\ \forall j\in\N{:} \qquad\mathbb{E}\bigl[e^{{D_j(\bx)}/{(C_4 T)}}|\Fill\bigr]\leq e^{j}.
\]
\end{lemma}

If $|V|\geq300$ Markov's inequality allows one to deduce that, for yet
another universal constant $L\equiv C_3 C_4$,
\[
\mathbb{P}\bigl(D_j(\bx)>t|\Fill\bigr)\leq e^{j - {t}/{(L \fT_{\IP(2,G)}(1/4))}}.
\]
Plugging this into (\ref{eqgettingthere112}) and Lemma~\ref
{lemchameleonmixing}, we obtain
\[
\dtv({\sL}[\bx^I_t]),\Unif((V)_k))\leq C_1 |V|^{3/2} e^{-c_1j} +
10|V| e^{j - {t}/{(L \fT_{\IP(2,G)}(1/4))}}.
\]
Since this inequality holds for all $j$, we can take
\[
j =\biggl\lfloor\frac{t}{2L \fT_{\IP(2,G)}(1/4)}\biggr\rfloor
\]
and obtain
\[
\dtv({\sL}[\bx^I_t],\Unif((V)_k))\leq K_0 |V|^{3/2} e^{-
{t}/{(2L \fT
_{\IP(2,G)}(1/4))}}
\]
with $K_0>0$ universal. Comparing with the definition of mixing time in
(\ref{eqdefmixingtime}) and noting that $\Unif((V)_k)$ is stationary for
$\IP(k,G)$ finishes the proof in the case $|V|\geq300$.

The case $|V|<300$---that is, $|V|$ bounded by a universal constant---can
be dealt with in several ways. For example, one may use the result
of Caputo et al.~\cite{CaputoEtAlInterchange} for the spectral gap of
$\IP(k,G)$ together with the standard \textit{lower} bound for $\fT_{\RW
(G)}(1/4)$ in terms of the spectral gap and the usual \textit{upper bound}
for $\fT_{\IP(k,G)}(\eps)$ in terms of its spectral gap (see, e.g.,
\cite{MontenegroTetaliMixingBook} for these standard bounds).
Alternatively, one may use Aldous and Fill's analysis (see Remark~\ref
{remaldousfill}) together with the inequality
\[
\mathbb{P}\bigl(M(\bx)>2i\fT_{\RW(G)}(1/4)\bigr)\leq\biggl(1-\frac
{1}{4|V|}\biggr)^i,
\]
which one can prove via Proposition~\ref{propmixingseparation} and a
few simple
calculations.
\end{pf}

\subsection{Outline of the missing steps}\label{secremainder}
We now summarize the main steps left in the proof.

\begin{enumerate}[(1)]
\item[(1)] In Section~\ref{secmiscellany} we collect several facts
about the
quantity $\ink$. The proof of Proposition~\ref{propdiscreteinkdecay}
on the decay
of $\mathbb{E}[1-\hatink^\bx_\ell/(|V|-k+1)\vert\Fill
]$, presented in the
\hyperref[app]{Appendix} (see Section~\ref{secproofdiscreteinkdecay}), relies on
results from
this section.
\item[(2)] Section~\ref{secinkconvergence} contains the proof of
Lemma~\ref{lemchameleonmixing}, which is based on an auxiliary result
on conditional
distributions (Lemma~\ref{lemchameleonmarginals}).
\item[(3)] Section~\ref{secdepinkingsarefast} bounds the right tail
of the
first depinking time in a chameleon process, and then uses this to
bound the exponential moment of the $j$th depinking time. This leads to
the key Lemma~\ref{lemexpmomentDj}, proven in Section~\ref
{secproofexpmomentDj}.
\end{enumerate}

%

\section{A miscellany of facts on ink}\label{secmiscellany}

In this section we prove several facts we will need about the quantity
$\ink_t$ introduced in (\ref{eqdeftotalink}). We will use the same
notation introduced in the proof of Lemma~\ref{lemkto2} (cf.
Section~\ref{secproofkto2}):
\begin{enumerate}[(1)]
\item[(1)] $\bx\in(V)_k$ is some fixed state;
\item[(2)] $(\bz,R,P,W)=(\bz,\{x\},\varnothing,V\setminus(\bO(\bz
)\cup\{
x\}))\in\sC_k(V)$ is the initial state corresponding to $\bx$ in the
sense of Proposition~\ref{propchameleonisrelated};
\item[(3)] $\ink_t^\bx$ is the total amount of ink in
$(z^C_t,R^C_t,P^C_t,W^C_t)$ (with the above initial state);
\item[(4)] $D_j(\bx)$ is the $j$th depinking time for this process;
\item[(5)] finally, $\hatink^\bx_j\equiv\ink^\bx_{D_j(\bx)}$.
\end{enumerate}
We will mostly omit $\bx$ from the notation in what follows. Most
proofs in this section follow by inspection, so we will be quite brief.

In principle the total number of number of depinking times could be
finite. We begin by showing that this is not the case.

\begin{proposition}\label{propinfmanydepinkings}The number of
depinking times is almost surely infinite.
\end{proposition}
\begin{remark}Notice that this only works because our \textit{definition} of a depinking time allows for ``trivial depinking times''
where there are only red and black (or only white and black) particles
left. This was noted in Box~\ref{boxchameleon}.
\end{remark}
\begin{pf*}{Proof of Proposition~\protect\ref{propinfmanydepinkings}} We
use the
following simple fact (proof omitted): there exists some $\delta>0$
such that each color-changing phase that starts with $\min\{
|R^C_t|,|W^C_t|\}>0$ will have a pinkening with probability $\geq
\delta
$, \textit{regardless of the past}. This implies that, given $s\geq0$,
the values of $|P^C_t|$ for $t\in[s,+\infty)$ will have a strictly
positive probability of increasing by the end of each color-changing
phase, at least until $|P^C_t|\geq\min\{|R^C_t|,|W^C_t|\}$. Since
$|P^C_t|$ can only decrease at depinking steps, this shows that
$|P^C_t|$ must continue to increase until $|P^C_t|\geq\min\{
|R^C_t|,|W^C_t|\}$, and the next time of the $2iT$ will be a depinking time.\vadjust{\goodbreak}
\end{pf*}

The next result follows by inspection.
\begin{proposition}[(Proof omitted)]\label{propinkbounded}$0\leq
\hatink
_{j}\leq|V|-k+1$ for all $j\in\N$, a.s.
\end{proposition}
%

We now compute the amount of change of ink in each step.
\begin{proposition}\label{propinkincrement}For $j\in\N$, $\hatink
_{j+1}\in\{\hatink_{j} + \Delta(\hatink_j),\hatink_j-\Delta
(\hatink
_j)\}$ a.s., where
%
%
\begin{equation}\label{eqdefDelta}\Delta(r)\equiv\biggl\lceil\frac
{\min
\{r,|V|-k+1-r\}}{3}\biggr\rceil \qquad  (r\in\N).
\end{equation}
Moreover, conditionally on $\{\hatink_\ell\}_{\ell=0}^j$, each
possibility is equally likely.
\end{proposition}
\begin{pf}Box~\ref{boxchameleon} shows that there are no pink
particles left in the system after depinking is performed. This implies
that $\hatink_j= \ink_{D_j} = |R^C_{D_j}|$. Moreover, since the total
number of nonblack particles is $|V|-k+1$, there must be $\hatink_j$
red and $|V|-k+1-\hatink_j$ white particles at time $D_j$.

A pinkening step decreases the number of red and white particles by $1$
each and increases the number of pink particles by $2$. However, no
pinkenings are performed if the number of pink particles is at least
the number of red or the number of white particles. In other words, the
number of pinkening steps until the next depinking is precisely the
\textit{smallest} $p$ satisfying $2p\geq\hatink_j-p$ or $2p\geq
|V|-k-1-\hatink_j-p$, which is $p=\Delta(\hatink_j)$ for $\Delta$
defined in (\ref{eqdefDelta}).

At the depinking step, the pink particles either all become white, or
they all become red. These possibilities corresponds to $\hatink
_{j+1}=\hatink_j-\Delta(\hatink_j)$ or $\hatink_{j+1}=\hatink
_j+\Delta
(\hatink_j)$, respectively. Which possibility will actually occur
depends on the value of the fair coin $d_i$, that is, flipped at the
depinking time $2iT=D_{j+1}$. It is easy to see that the coin is
independent of $\{\hatink_{\ell}\}_{\ell\leq j}$, and this implies that
both possibilities are equally likely.
\end{pf}

The next lemma summarizes the above sequence of propositions and adds a
useful remark.

\begin{lemma}\label{lemsummaryhatink}The sequence $\{\hatink_j\}
_{j\geq0}$ is a Markov chain with initial state $\hatink_0=1$,
absorbing states at $0$ and $|V|-k+1$ and transition probabilities
given by
%
%
\begin{eqnarray}\label{eqsummarytransitionprobs}p(a,b)\equiv\frac
{1}{2}\bigl({\mathbb{I}_{\{b=a+\Delta(a)\}}} + {\mathbb{I}_{\{
b=a-\Delta(a)\} }}\bigr)
\nonumber
\\[-8pt]
\\[-8pt]
\eqntext{(a,b\in\{0,1,2,\ldots,|V|-k+1\}).}
\end{eqnarray}
Moreover, it is almost surely absorbed in finite time in either $0$ or
$|V|-k+1$. Finally, the event
%
%
\begin{equation}\label{eqdefFill}\Fill\equiv\Bigl\{\lim_{j\to
+\infty
}\hatink_{j} = |V|-k+1\Bigr\}
\end{equation}
has probability $1/(|V|-k+1)$.
\end{lemma}
\begin{remark}\label{remFill}The event $\Fill$ corresponds to the
number of red particles converging to $|V|-k+1$, that is, that there
are only black and red particles at all large enough times, or,
equivalently, to red ink filling up all available space. Notice that we
can rewrite
\[
\Fill\equiv\Bigl\{\lim_{t\to+\infty}\ink_t=|V|-k+1\Bigr\},
\]
which is the form that appears in the proof of Lemma~\ref{lemmixfor2}.
\end{remark}

\begin{pf*}{Proof of Lemma~\protect\ref{lemsummaryhatink}} The first
sentence is
obvious given the sequence of propositions; only notice that $p(a,a)=1$
if $a\in\{0,|V|-k+1\}$. We omit the trivial proof of the next
assertion, which implies $\hatink_{\infty}\equiv\lim_{j\to+\infty
}\hatink_{j}\in\{0,|V|-k+1\}$.

Now notice that the increments of $\hatink_j$ are unbiased; that
implies that this process is also a martingale. We thus have
\begin{eqnarray*}
\mathbb{P}(\Fill) & = &\mathbb{P}(\hatink_{\infty
} = |V|-k+1)\\
\mbox{($\hatink_{\infty}\in\{0,|V|-k+1\}$)} &=& \frac{\mathbb
{E}[\hatink_{\infty}]}{|V|-k+1}\\
\mbox{($\{\hatink_j\}_{j\in\N}$ bounded, cf. Prop. \ref
{propinkbounded})} &=& \frac{\lim_{j\to+\infty}\mathbb{E}
[\hatink_{j}]}{|V|-k+1}\\
\mbox{(mart. property + $\hatink_0=1$)} &=& \frac{\mathbb{E}
[\hatink_0]}{|V|-k+1} = \frac{1}{|V|-k+1}.
\end{eqnarray*}
\upqed\end{pf*}

We will need one final lemma before we proceed.

\begin{lemma}\label{lemztindepFill}For all $\bb\in(V)_{k-1}$ and
$t\geq0$,
\[
\mathbb{P}(\{\bz^C_t=\bb\}\cap\Fill) = \frac{\mathbb{P}(\bz
^C_t=\bb)}{|V|-k+1}.
\]
\end{lemma}
\begin{pf}This follows from the previous lemma if we can show that
$\Fill$ and $\bz^C_t$ are independent. To see this, simply notice that
$\Fill$ is entirely determined by the coin flips $d_i$ performed at the
depinking times, whereas the value of $\bz^C_t$ does not at all depend
on these coin flips.
\end{pf}

\begin{remark}\label{remproofdiscreteinkdecay}It transpires from
the above that the chameleon process conditioned on $\Fill$ is the same
as the unconditional process, \textit{except that} the coin flips $d_i$
performed at depinking times are biased. This remark will be useful in
the proof of Lemma~\ref{lemexpmomentDj} in Section~\ref
{secproofexpmomentDj}.
\end{remark}

\section{Convergence to stationarity in terms of ink}\label
{secinkconvergence}

In this section we will prove Lemma~\ref{lemchameleonmixing}, used in the
proof of Lemma~\ref{lemkto2} (cf. Section~\ref{secproofkto2}), in
which we show that
the amount of ink in the system can be used to bound the distance to
the stationary distribution. We start with a preliminary result on marginals.

\subsection{The convergence to equilibrium of
conditional distributions}\label{secmarginals}

We will again use Notational convention~\ref{notlastcoordinate},
whereby any $\bx=(\bx
(1),\ldots,\bx(k))\in(V)_k$ is written as a pair $\bx=(\bz,x)$ with
$\bz
=(\bx(1),\ldots,\bx(k-1))$ and $x=\bx(k)$.

Let $\bx=(\bz,x)\in(V)_k$ and consider the $\IP(k,G)$ process $\{
\bx
^I_t\}_{t\geq0}$. Set $R=\{x\}$, $P=\varnothing$ and $W=V\setminus
(\mathbf{O}(\bz)\cup\{x\})$ and recall from Proposition~\ref
{propchameleonisrelated}
that the chameleon process $\{(\bz^C_t,R^C_t,P^C_t,W^C_t)\}_{t\geq0}$
satisfies
%
%
\begin{equation}\label{eqchameleonisrelated2}\quad\forall t\geq0,
\forall
\bb=(\bc,b)\in(V)_k\qquad  \mathbb{P}(\bx^I_t = \bb) = \mathbb
{E}
\bigl[{\mathbb{I}_{\{\bz^C_t = \bc\} }}\ink^{\bx}_t(b)\bigr],
\end{equation}
where (as before) we use $\ink^{\bx}_t(\cdot)$ to denote the amount
of ink in this chameleon process corresponding to $\bx$. The
following lemma relates the \textit{total} amount of ink in this
process to the near-uniformity of $x^I_t$ conditionally on
$\bz^I_t$.

\begin{lemma}\label{lemchameleonmarginals}Given $\bx=(\bz,x)\in
(V)_k$, let $\tilde{\bx}^I_t=(\bz^I_t,\tilde{x}^I_t)$ where,
conditionally on $\bz^I_t$, $\tilde{x}^I_t$ is uniform over
$V\setminus\mathbf{O}(\bz^I_t)$. Then
\[
\dtv({\sL}[\bx^I_t]),{\sL}[\tilde{\bx}^I_t])\leq\mathbb{E}
\biggl[1 - \frac{\ink^{\bx}_t}{|V|-k+1} \Big| \Fill\biggr]
\]
where $\Fill$ is the event defined in Lemma~\ref{lemsummaryhatink}
(see also
Remark~\ref{remFill}).
\end{lemma}

\begin{pf}We have seen that $\bz^I_t$ and $\bz^C_t$ have the same
distribution; cf. the proof of Proposition~\ref
{propchameleonisrelated}. We deduce that
\begin{eqnarray*}
\forall t\geq0, \forall\bb=(\bc,b)\in(V)_k \qquad \mathbb{P}(\bz
^I_t=\bc ,\tilde {x}^I_t=b) &=& \frac{\mathbb{P}(\bz^C_t=\bc
)}{|V|-k+1}\\
&=&\mathbb{P}(\{\bz^C_t=\bc\} \cap \Fill),
\end{eqnarray*}
where the last equality follows from Lemma~\ref{lemztindepFill}. On
the other
hand, (\ref{eqchameleonisrelated2}) implies
%
%
\begin{equation}\label{eqchameleonisrelated3}\qquad\forall t\geq0,
\forall
\bb=(\bc,b)\in(V)_k \qquad \mathbb{P}(\bx^I_t = \bb) \geq\mathbb
{E}
\bigl[{\mathbb{I}_{\{\bz^C_t = \bc\}\cap\Fill}}\ink^{\bx}_t(b)\bigr].
\end{equation}
We deduce that
\begin{eqnarray*}
&&\forall t\geq0, \forall\bb=(\bc,b)\in(V)_k\\[-2pt]
&&\qquad\bigl(\mathbb{P}(\bz^I_t=\bc,\tilde{x}^I_t=b)-\mathbb{P}(\bz^I_t=\bc
,x^I_t=b)\bigr)_+\\[-2pt]
&&\qquad\qquad\leq
\bigl(\mathbb{E}\bigl[{\mathbb{I}_{\{\bz^C_t = \bc\}\cap\Fill
}}\bigl(1-\ink
^{\bx}_t(b)\bigr)\bigr]
\bigr)_+ \\[-2pt]
&&\qquad\qquad
= \mathbb{E}\bigl[{\mathbb{I}_{\{\bz^C_t = \bc\}\cap
\Fill }}\bigl(1-\ink^{\bx}_t(b)\bigr)\bigr]\qquad\mbox{since the integrand is
$\geq0$.}
\end{eqnarray*}

We now combine this with formula (\ref{eqdtvsumpositivepart}) for
$\dtv
(\cdot,\cdot\cdot)$.
\begin{eqnarray*}
\dtv({\sL}[\bx^I_t],{\sL}[\tilde{\bx}^I_t])
&\leq&\sum
_{\bb=(\bc,b)\in(V)_k}\mathbb{E}\bigl[{\mathbb{I}_{\{\bz^C_t =
\bc\} \cap\Fill}}\bigl(1-\ink^{\bx}_t(b)\bigr)\bigr]\\[-2pt]
&=& \sum_{\bc\in
(V)_{k-1}}\mathbb{E}\biggl[{\mathbb{I}_{\{\bz^C_t = \bc\}\cap\Fill
}}
\sum_{b\in V\setminus\mathbf{O}(\bc)}\bigl(1-\ink^{\bx}_t(b)\bigr)\biggr]\\[-2pt]
\mbox{[sum over \textit{b} $+$ (\ref{eqtotalink2})]} &=& \sum_{\bc\in
(V)_{k-1}}\mathbb{E}\bigl[{\mathbb{I}_{\{\bz^C_t = \bc\}\cap
\Fill }}(|V|-k+1-\ink^{\bx}_t)\bigr]\\[-2pt]
\mbox{(sum over $\bc$)} &=&
\mathbb{E}[{\mathbb{I}_{\Fill}}(|V|-k+1-\ink^{\bx}_t)
]\\[-2pt]
\mbox{(apply Lemma~\ref{lemsummaryhatink})} &=& 1 - \frac{\mathbb
{E}
[\ink^{\bx}_t\vert\Fill]}{|V|-k+1}.
\end{eqnarray*}
\upqed\end{pf}

\subsection{Distance to the stationary distribution in terms of
ink}\label{secproofchameleonmixing}
\mbox{}
\begin{pf*}{Proof of Lemma~\protect\ref{lemchameleonmixing}} We will prove the
following stronger inequality:
%
%
\setcounter{equation}{0}
\begin{equation}\label{eqchameleonarbitraryinitial}\sup_{\bx,\by
\in
(V)_k}\dtv({\sL}[\bx_t^I],{\sL}[\by_t^I])\leq
2k\sup_{\bw\in(V)_k}\mathbb{E}\biggl[ 1 - \frac{\ink^{\bw
}_t}{|V|-k+1} \Big| \Fill\biggr],
\end{equation}
which
implies the lemma by convexity.

Declare two states $\bu,\bv\in(V)_k$ to be \textit{adjacent}
($\bu\sim\bv$) if they differ at precisely one coordinate: that is, there
exists an $i\in[k]$ with $\bu(i)\neq\bv(i)$ and $\bu(r)=\bv(r)$ for
$r\in[k]\setminus\{i\}$. We first bound
$\dtv({\sL}[\bx_t^I],{\sL}[\by_t^I])$ for \textit{adjacent} $\bx\sim
\by$.

One can assume without loss of generality that $\bx$ and $\by$
differ precisely at the $k$th coordinate. Using the notation from
Section~\ref{secchameleonisrelated}, we write $\bx=(\bz,x)$
and $\by=(\bz,y)$ for $\bz\in(V)_{k-1}$ and $x\in V\setminus\mathbf{O}(\bz)$. Defining $\tilde{\bx}_t^I=(\bz^I_t,\tilde{x}^I_t)$ as in
Section~\ref{secmarginals} and $\tilde{\by}^I_t$ similarly, we see that
${\sL}[\tilde{\bx}_t^I]={\sL}[\tilde{\by}_t^I]$ for all $t\geq
0$. We
deduce that
%
%
\begin{eqnarray}
\dtv({\sL}[\bx_t^I],{\sL}[\by_t^I])
&\leq
&\dtv({\sL}[\bx_t^I],{\sL}[\tilde{\bx}_t^I]) +
\dtv({\sL}[\by_t^I],{\sL}[\tilde{\by}_t^I]) \nonumber\\
&&{} +
\dtv({\sL}[\tilde{\bx}_t^I],{\sL}[\tilde{\by}_t^I])\nonumber\\
\mbox{(3rd. term${}=0$)}&=&
\dtv({\sL}[\bx_t^I],{\sL}[\tilde{\bx}_t^I]) +
\dtv({\sL}[\by_t^I],{\sL}[\tilde{\by}_t^I])
\nonumber
\\[-8pt]
\\[-8pt]
\nonumber
\mbox{(use Lemma~\ref{lemchameleonmarginals})} &\leq& \mathbb{E}\biggl[
1 -
\frac{\ink^{\bx}_t}{|V|-k+1} \Big| \Fill\biggr]\\
&&{} + \mathbb
{E}\biggl[ 1 - \frac{\ink^{\by}_t}{|V|-k+1} \Big| \Fill
\biggr]\nonumber\\
\label{eqchameleonadjacent} &\leq& 2\sup_{\bw\in(V)_k}\mathbb
{E}\biggl[ 1 - \frac{\ink^{\bw}_t}{|V|-k+1} \Big| \Fill
\biggr].\nonumber
\end{eqnarray}

Now consider $\bx,\by\in(V)_k$ arbitrary. One can find a sequence
$\{\bx[i]\}_{i=0}^r\subset(V)_k$ with $r\leq2k$ and
\[
\bx[0]=\bx\sim\bx[1]\sim\bx[2]\sim\cdots\sim\bx[r]=\by.
\]
The triangle
inequality gives
\[
\dtv({\sL}[\bx_t^I],{\sL}[\by_t^I]) =
\dtv({\sL}[\bx[0]_t^I],{\sL}[\bx[r]_t^I])\leq\sum_{i=1}^r
\dtv({\sL}[\bx[i-1]_t^I],{\sL}[\bx[i]_t^I]).
\]
Applying
(\ref{eqchameleonadjacent}) to each adjacent pair $\bx[i-1],\bx[i]$
gives (\ref{eqchameleonarbitraryinitial}).
\end{pf*}

\section{Depinkings are fast}\label{secdepinkingsarefast}

The results in this section lead to the key Lem\-ma~\ref
{lemexpmomentDj}. We
first show that, in the first two phases of
the chameleon process---a constant color and a color-changing phase---,
the number of red particles decreases in
expectation by a constant factor.

\begin{lemma}[(Proven in Section~\protect\ref{secprooffirsttworounds})]\label
{lemfirsttworounds}Consider a modified chameleon process where one
drops condition (2) for a pinkening step; cf. Box~\ref{boxchameleon}.
Assume also that $k\leq|V|/2$, $|V|\geq300$ and that the initial
state $(\bz,R,P,W)\in\sC_k(V)$ with $|P|<|R|\leq|W|$. If the phase
length parameter $T$ satisfies
\[
T\geq20\fT_{\IP(2,G)}(1/4),
\]
then
\[
\mathbb{E}[|R^C_{2T_-}|]\leq(1-c)|R|,
\]
where $c=1/1000>0$.
\end{lemma}

With this, we will show that the first depinking time has
an exponential moment.

\begin{lemma}[(Proven in Section~\ref{secprooffirstdepinking})]\label
{lemfirstdepinking}Consider a chameleon process (without the
modification in the previous lemma) with $|V|\geq300$ and $k\leq
|V|/2$, started from an initial state $(\bz,R,P,W)\in\sC_k(V)$ with
$|P|=\varnothing$. There exists a universal constant $K>0$ such that if
the phase length parameter $T$ satisfies
$T\geq20 \fT_{\IP(2,G)}(1/4)$, the first depinking
time $D_1$ of this process satisfies
\[
\mathbb{E}\bigl[e^{{D_1}/{(K T)}}\bigr]\leq e.
\]
\end{lemma}

In Section~\ref{secproofexpmomentDj} we deduce Lemma~\ref
{lemexpmomentDj} from
Lemma~\ref{lemfirstdepinking}.

\subsection{Loss of red particles in the two first
phases}\label{secprooffirsttworounds}
\mbox{}
\begin{pf*}{Proof of Lemma~\ref{lemfirsttworounds}} Note that
there is no depinking at time $t=0$, since
there are less pink particles than white or red ones in the state $(\bz
,R,P,W)$. Finally, the conditions on $P,W,R$ and $k$ imply
\[
3|W|\geq|R|+|P|+|W| = |V|-k+1\geq|V|/2\quad\Rightarrow\quad|W|\geq|V|/6.
\]

The interval $(0,T]$ is a constant-color phase where
black, red and white particles are simply moved around. Lemma~\ref
{lemjustlikeI} shows that the state of the process at time~$T$ is
given by
\[
(\bz^C_T,R^C_T,P^C_T,W^C_T) = (I(\bz),I(R),I(P),I(W)),
\]
where $I=I_{(0,T]}=I_T$ is the map obtained from the modified
chameleon construction in Section~\ref{secchameleon}. We will need the
following properties later on:
\begin{proposition}[(Proven in Section~\ref{secproofI})]\label
{propI}For all
$(a,b)\in(V)_2$ and $S,L\subset V$ with $S\times L\subset(V)_2$,
$|L|\geq|V|/12$,
\begin{eqnarray*}
\mathbb{P}\bigl((a,b)\in I(S)\times I(S)\bigr)&\leq&\mathbb{P}(a\in I(S))
\biggl(\frac{|S|}{|V|}+2^{-10}\biggr).
\\
\mathbb{P}\bigl((a,b)\in I(S)\times I(L)\bigr)&\geq&\frac{|S||L|}{|V|^2}
(1-2^{-9})\geq\frac{|S|}{13|V|}.
\end{eqnarray*}
\end{proposition}
\begin{remark}\label{remexplainconstantcolor}The intuitive meaning
of this is that $(R^C_T,W^C_T)$
are close to uniform in terms of correlations of ``pairs of
particles'' at the end of the constant-color phase, and this will
only hold because $T=\Omega(\fT_{\IP(2,G)})$. Morris's original
argument for $(\Z/L\Z)^d$ could instead rely on good estimates for
transition probabilities for single-particle random walks. We note that
we need the negative correlation property in the proof of this proposition.
\end{remark}

In the time interval $(T,2T)$, each time $T<\tau_m<2T$ may or may not
be a
pinkening time, depending on whether pinkening condition $(1)$ is
satisfied. We will nevertheless consider the maps
%
%
\begin{equation}\label{eqdeftildeI}\tilde{I}_t\equiv I_{(T,t]},\qquad
T\leq t<2T \mbox{; cf. the definition in Section~\ref{secmodgraphical}}.
\end{equation}
We emphasize that $\tilde{I}_t$ does \textit{not}
correspond directly to the evolution of the chameleon process in the
time interval $(T,2T)$. Propositions~\ref{propIndependent} and \ref
{propmodisthesame} imply:
\begin{proposition}[(Proof omitted)]\label{proptildeI}$\{\tilde{I}_t\}
_{T<t<2T}$ is independent from $I$, and so are all the points of the
Poisson process $\{\tau_n\}_n$ in the interval $(T,2T)$ and all
markings $e_n,c_n$ corresponding to these points.
\end{proposition}

We need a new definition before we proceed. Let $a\in V$ be given.
Let $\phi_a$ be the \textit{first time} of the form\vadjust{\goodbreak} $\tau_m$ with
$T<\tau_m\leq2T$ for which $a\in e_m$; if no such time exists, let
$\phi_a=+\infty$. If $\phi_a<+\infty$, there exists a vertex $b\in
V$ such that the edge $e_m$ just mentioned has
$a=\tilde{I}_{{\phi_a}_-}(a)$ and $\tilde{I}_{{\phi_a}_-}(b)$ as
endpoints immediately prior to time $\phi_a$. We set $F_a\equiv b$ in
that case, or $F_a\equiv*$ if
$\phi_a=+\infty$. The following simple claim is essential to what
follows.

\begin{claim}\label{claimnoofpinkening}The number of pinkening steps
performed in time interval $(T,2T)$ is \textit{at least} the number of
$b\in I(W)$ such that $F_a=b$ for some \mbox{$a\!\in\! I(R)$}.\looseness=-1
\end{claim}
\begin{pf}
Let $b\in I(W)$. Given the rules for
color-changing phases (cf. Box~\ref{boxchameleon}), the particle at
that location will move in the time interval $(T,2T)$ according to
$\tilde{I}_t$ until the first time $t=\tau_m\in(T,2T)$ such that
$\tilde{I}_{t_-}(b)\in e_m$ and the other endpoint of $e_m$ is white
(if such a time exists). Now if $a\in I(R)$ satisfies $F_a=b$ and $\tau
_m=\phi_a$, we have $\{\tilde{I}_{t_-}(b),a\}=e_m$ and $a$ must still
be red at time $(\phi_a)_-$, since it was not contained in an edge
before in this phase. It follows that the particle started from $b$
must become pink by time $\phi_a$.
\end{pf}

The claim implies
%
%
\begin{eqnarray}\label{eqtowardstheendRC1}
|R^C_{2T_-}| &=& |R| -
\mbox{$\#$ of pinkening steps in }(T,2T] \\[-2pt]
\label{eqtowardstheendRC11}
&\leq& |R| - \sum_{b\in I(W)} {\mathbb
{I}_{\bigcup _{a\in I(R)}\{F_a = b\}}}.
\end{eqnarray}
The sum in the RHS satisfies
%
%
\begin{eqnarray}\label{eqZopa}\sum_{b\in I(W)} {\mathbb{I}_{\bigcup
_{a\in I(R)}\{ F_a = b\}}}&\geq&\sum_{a\in I(R),b\in I(W)}{\mathbb
{I}_{\{F_a=b\}}}
\nonumber
\\[-9pt]
\\[-9pt]
\nonumber
&&{}- \sum_{\{a,a'\}
\subset I(R),b\in I(W)}{\mathbb{I}_{\{F_a=b\}}} {\mathbb{I}_{\{
F_{a'}=b\}}},
\end{eqnarray}
and we obtain
\begin{eqnarray}\label{eqtowardstheendRC2}
&&\mathbb{E}
[|R^C_{2T_-}|-|R|]\nonumber\\[-2pt]
&&\qquad\leq
-\sum_{(a,b)\in(V)_2}\mathbb{P}\bigl(a\in I(R),b\in I(W),F_a=b\bigr)
\\[-2pt]
&&\qquad\qquad{} + \sum_{\{a,a'\}\subset V,b\in V}\mathbb{P}\bigl(a,a'\in I(R),b\in
I(W),F_a=b,F_{a'}=b\bigr).\nonumber
\end{eqnarray}

The event $\{F_a=b\}$ is entirely determined by the points of the
marked Poisson process and by the coin flips performed in the time
interval $(T,2T)$, and therefore
is independent of $I$; cf. Proposition~\ref{proptildeI}. We deduce
\begin{eqnarray}\label{eqtowardstheendRC3}
&&\sum_{(a,b)\in
(V)_2}\mathbb{P}\bigl(a\in I(R),b\in I(W),F_a=b\bigr) \nonumber\\[-2pt]
&&\qquad\hspace*{76pt}= \sum_{(a,b)\in
(V)_2}\mathbb{P}\bigl(a\in I(R),b\in I(W)\bigr)\mathbb{P}(F_a=b)
\nonumber
\\[-8pt]
\\[-8pt]
\nonumber
&&\mbox{(use Proposition~\ref{propI})
} \geq\frac{|R|}{13|V|} \sum_{(a,b)\in(V)_2}\mathbb{P}(F_a=b) \\
&&\qquad\hspace*{76pt}=
\frac{|R|}{13|V|} \sum_{a\in V}\mathbb{P}(F_a\neq*).\nonumber
\end{eqnarray}
For a
given $a\in V$, $\mathbb{P}(F_a=*)$ is the probability that there is no
$T<\tau_n<2T$ with $e_n\ni a$. Notice that this is at most the
probability that $\tilde{I}_{2T}(a)=a$: $a$~cannot move if there is
no edge $e_n\ni a$ with $T<\tau_n\leq2T$. We deduce
\[
\mathbb{P}(F_a\neq*)\geq1 - \mathbb{P}\bigl(\tilde{I}_{2T}(a)=a\bigr) =
\mathbb{P}(a^R_{T}\neq a),
\]
where $\{a^R_t\}_{t\geq0}$ is a realization of $\RW(G)$ started from
$a$. By the contraction principle and Proposition~\ref{propmixingpowersof2},
\[
T\geq20\fT_{\IP(2,G)}(1/4)\geq20\fT_{\RW(G)}(1/4)\geq\fT_{\RW
(G)}(2^{-20}),
\]
which implies
\[
\mathbb{P}(a^R_{T}\neq a)\geq1 - \frac{1}{|V|} - 2^{-20}\geq\frac
{13}{14}\qquad \mbox{since $|V|\geq300$}.
\]
We deduce from (\ref{eqtowardstheendRC3}) that
%
%
\begin{equation}\label{eqteRC21!}\sum_{(a,b)\in V^2}\mathbb
{P}\bigl(a\in I(R),b\in I(W),F_a=b\bigr) \geq\frac{|R|}{14}.
\end{equation}

We now consider the second sum in the RHS of
(\ref{eqtowardstheendRC2}). As before, we notice that $\{
F_a=b,F_{a'}=b\}
$ is independent of $I$ and therefore
\begin{eqnarray*}
&&\sum_{\{a,a'\}\subset V,b\in V\setminus\{a,a'\}
}\mathbb{P}\bigl(a,a'\in I(R),b\in I(W),F_a=b,F_{a'}=b\bigr) \\
&&\qquad= \sum_{\{
a,a'\}\subset V,b\in
V\setminus\{a,a'\}}\mathbb{P}\bigl(a,a'\in I(R),b\in I(W)\bigr)\mathbb
{P}(F_a=b,F_{a'}=b)\\
&&\qquad\leq\sum_{\{a,a'\}\subset V,b\in
V\setminus\{a,a'\}}\mathbb{P}\bigl(a,a'\in I(R)\bigr)\mathbb{P}(F_a=b,F_{a'}=b).
\end{eqnarray*}

We \textit{claim} that:
\begin{claim}[(Proven in Section~\ref
{secproofgraphicalargument})]\label
{claimgraphicalargument}For all $(a,a',b)\in(V)_3$,
\[
\mathbb{P}(F_a=b,F_{a'}=b)\leq\mathbb{P}(F_a=a',F_{a'}=b) +
\mathbb{P}(F_{a'}=a,F_a=b).
\]
\end{claim}
\ignore{
\begin{remark}Claim~\ref{claimgraphicalargument} will help us control the
probability of two red particles trying to pinken when ``bumping into''
the same white particle. The next equation will show that the
probability that $F_a=F_{a'}$ is bounded by $\mathbb
{P}(F_a=a')+\mathbb{P}(F_{a'}=a)$. We were not able to prove a similar
result for \textit{general graphs $G$} without rule $\# 3$ for a {\sc pinkening time} (cf.
Boxes~\ref{boxpinkening} and~\ref{boxnewpinkening}): that is, the
fact that pinkening can only happen when a vertex is \textit{first}
touched by an edge is what allows us to show that $F_a=F_{a'}$ is not
too likely. By contrast, if \textit{any contact} between red and white
could result in pinkening, problems would occur eg. in a graph with a
white vextex of very large degree with many red neighbors.
\end{remark}
}

Summing up over $b$ above gives at most $\mathbb{P}(F_a=a')+\mathbb
{P}(F_{a'}=a)$ in
the RHS. Therefore,
\[
\hspace*{-130pt}\sum_{\{a,a'\}\subset V,b\in V\setminus\{a,a'\}
}\mathbb{P}\bigl(a,a'\in I(R)\bigr)\mathbb{P}(F_a=b,F_{a'}=b)\vspace*{-9pt}
\]
\begin{eqnarray*}
&\leq&\sum_{\{
a,a'\}\subset V}\mathbb{P}\bigl(a,a'\in I(R)\bigr)
\bigl(\mathbb{P}(F_a=a')+\mathbb{P}(F_{a'}=a)\bigr)\\
&=&\sum_{(a,a')\in(V)_2}\mathbb{P}\bigl(a\in I(R),a'\in I(R)\bigr)\mathbb
{P}(F_a=a') \\
\mbox{(apply Prop.~\ref{propI})}&=& \biggl(\frac{|R|}{|V|}+2^{-10}
\biggr)\sum
_{(a,a')\in
(V)_2}\mathbb{P}\bigl(a\in I(R)\bigr)\mathbb{P}(F_a=a')\\
\biggl(\bigcup_{a'}\{
F_a=a'\}=\{
F_a\neq
*\}\biggr)  &= &\biggl(\frac{|R|}{|V|}+2^{-10}\biggr)\sum_{a\in V}\mathbb
{P}\bigl(a\in I(R)\bigr)\mathbb{P}(F_a\neq*)\\
\bigl(\mathbb{P}(F_a\neq
*)\leq1\bigr) &\leq&
\biggl(\frac
{|R|}{|V|}+2^{-10}\biggr)\sum_{a\in V}\mathbb{P}\bigl(a\in I(R)\bigr)\\
& =&
\biggl(\frac
{|R|}{|V|}+2^{-10}\biggr)\mathbb{E}[|I(R)|] \\
& =&
\biggl(\frac
{|R|}{|V|}+2^{-10}\biggr)|R|\qquad\mbox{\hspace*{-3.5pt}since $I=I_{(0,T]}$ is a bijection.}
\end{eqnarray*}
Plugging this
equation and (\ref{eqteRC21!}) into (\ref{eqtowardstheendRC2}) we obtain
\begin{eqnarray}\label{eqtowardstheendRC4}
\mathbb{E}[|R^C_{2T}|-|R|]&\leq&|R|
\biggl(\frac{|R|}{|V|}+2^{-10}-\frac{1}{14}\biggr)
\nonumber
\\[-8pt]
\\[-8pt]
\nonumber
 &\leq&
-|R|/30\qquad\mbox{if }|R|\leq|V|/28.
\end{eqnarray}
If $|R|>|V|/28$, we can
still find a subset $R_0\subset R$ of size $|R_0|=\lfloor|V|/28\rfloor
$. Since
%
%
\begin{equation}\label{eqZopa2}\sum_{b\in I(W)} {\mathbb{I}_{\{
\exists a\in I(R) :  F_a = b\}}}\geq\sum_{b\in I(W)}{\mathbb
{I}_{\{\exists a\in I(R_0) : F_a=b\}}},
\end{equation}
we may repeat the reasoning presented from (\ref{eqZopa}) onwards,
replacing $R$ by $R_0$, to deduce that
\[
\mathbb{E}[|R^C_{2T_-}|-|R|]\leq-\frac{|R_0|}{30}.
\]
We now note that, since $|V|\geq300$,
\[
|R_0|\geq\frac{|V|}{30} - 1 \geq\frac{3|V|}{100}\geq\frac{3|R|}{100}
\]
since $|R|\leq|V|$. We deduce that
\[
\mathbb{E}[|R^C_{2T_-}|-|R|]\leq- \frac{|R|}{1000}\qquad\mbox{if }|R|>|V|/28,
\]
which gives the lemma together with (\ref{eqtowardstheendRC4}).
\end{pf*}
%
\subsubsection{\texorpdfstring{Proof of the required estimates for the $I$ map (Proposition~\protect\ref{propI})}
{Proof of the required estimates for the $I$ map (Proposition 9.1)}}\label{secproofI}
\mbox{}
\begin{pf*}{Proof of Proposition~\protect\ref{propI}}Recall that $T\geq
20\fT_{\IP
(2,G)}(1/4)$, therefore $T\geq2\fT_{\IP(2,G)}(2^{-10})$ by
Proposition~\ref{propmixingpowersof2}. By the contraction principle
\cite
{AldousFillRWBook}, this also implies that $T\geq\fT_{\RW(G)}(2^{-10})$.

Recall that $I=I_{(0,T]}$ as in the construction of the modified
chameleon process. This implies that for any set $S$, $I(S)$ has the
law of $\EX(|S|,G)$ started from $S$. We deduce
\begin{eqnarray*}
\mathbb{P}\bigl((a,b)\in I(S)\times I(S)\bigr) &=& \mathbb
{P}(\{a,b\}\subset S^I_T)\\
\mbox{(negative correlation, Lemma~\ref{lemNEGCORR1})} &\leq& \mathbb{P}(a\in S_T^I)\mathbb{P}(b\in
S^I_T)\\
\mbox{($\sL[I]=\sL[I^{-1}]$, Proposition~\ref{propbackwards})} & = & \mathbb{P}\bigl(a\in I(S)\bigr)\mathbb{P}(b^I_T\in
S)\\
\bigl(T\geq \fT_{\RW(G)}(2^{-10})\bigr) &\leq& \mathbb{P}\bigl(a\in I(S)\bigr)\biggl(\frac
{|S|}{|V|}+2^{-10}\biggr).
\end{eqnarray*}
As for the other inequality in the proposition, we have
\begin{eqnarray*}
\mathbb{P}\bigl((a,b)\in I(S)\times I(L)\bigr) &=& \mathbb
{P}\bigl((a^I_T,b^I_T)\in S\times L\bigr)\\
\bigl(\mbox{take }\bx=(a,b)\bigr) &=&
\mathbb{P}(\bx^I_T\in S\times L)\\
\mbox{(*)} & \geq& (1-2^{-9})^2 \frac{|S\times
L|}{|(V)_2|}\\
&\geq& (1-2^{-8}) \frac{|S| |L|}{|V|^2},
\end{eqnarray*}
where $(*)$ follows from the symmetry of the transition rates of $\IP
(2,G)$, the fact that $T\geq2\fT_{\IP(2,G)}(2^{-10})$ and
Proposition~\ref{propmixingseparation}. We note that $|L|/|V|\geq
(1/12)$ and $1-2^{-8}\geq
12/13$ to finish the proof.
\end{pf*}

\subsubsection{\texorpdfstring{Proof of claim on $F_a$ (Claim~\protect\ref{claimgraphicalargument})}
{Proof of claim on $F_a$ (Claim 2)}}\label{secproofgraphicalargument}
\mbox{}
\begin{pf*}{Proof of Claim~\protect\ref{claimgraphicalargument}}
It suffices to show that for $(a,a',b)\in(V)_3$,
%
%
\begin{equation}\label{eqgraphicalargument1}\qquad\mathbb
{P}(F_a=b,F_{a'}=b,\phi _a\leq\phi_{a'}) =\mathbb
{P}(F_a=b,F_{a'}=a,\phi_a\leq\phi_{a'}).
\end{equation}
Let $L_b,R_b$ denote the events appearing
in the LHS and RHS of (\ref{eqgraphicalargument1}) (resp.). We
present a simple measure-preserving mapping $\Phi$, which acts on
\[
(\sP,\{e_n\}_{n},\{c_n\}_n,\{d_i\}_{i}),
\]
that maps $L_b$
into $R_b$ and vice-versa. We describe $\Phi$ in words: all values of
$d_i$, $T<\tau_j\leq2T$ and all corresponding $e_j$ and $c_j$, \textit{except} for the following modification: if $\tau_m=\phi_a$, we flip
the value of $c_m$ to $c'_m=1-c_m$.

Let us check that $\Phi$ has the desired properties. $\Phi$ is
clearly measure-preserving, since $\phi_a$ is a stopping time that
is independent of the value $c_m$ of the flipped coin.

Now suppose $\{\hat{I}_t\}_{T<t\leq2T}$ is defined precisely as
$\{\tilde{I}_t\}_{T<t\leq2T}$, but with $c_m$ flipped. It is easy to
see that $\phi_a,\phi_{a'}$ retain their
values and that the random variable $\hat{F}_a$ corresponding to $F_a$
in the $\hat{I}$ process satisfies $\hat{F}_a=F_a$. The two
processes coincide for any time $T<t<{\phi_a}$. If $L_b$ holds, we have
$\phi_a=\tau_j<2T$ for some $j$, and the endpoints of $e_j$ are $a$ and
$\tilde{I}_{{\phi_a}_-}(F_a) = \tilde{I}_{{\phi_a}_-}$ (by
definition of $F_a$). Since the coin flips used for
$\tilde{I}_{\phi_a}$ and $\hat{I}_{\phi_a}$ are opposite, we have
\[
(\hat{I}_{\phi_a}(a),\hat{I}_{\phi_a}(b)) = (\tilde{I}_{\phi
_a}(b),\tilde{I}_{\phi_a}(a))
\]
whereas $\tilde{I}_{\phi_a}(c)=\hat{I}_{\phi_a}(c)$ for all $c\in
V\setminus\{a,F_a\}$.
\[
\mbox{Under $L_b$, }\forall t\in[\phi_a,2T]{:}\qquad  (\hat
{I}_{t}(a),\hat
{I}_{t}(b)) = (\tilde{I}_{t}(b),\tilde{I}_{t}(a)).
\]
Under $L_b$, the edge $e_\ell$ corresponding to $\tau_\ell=\phi_{a'}$
was of
the form $e_\ell=\{a',\break \tilde{I}_{{\phi_{a'}}_-}(b)\}$. This implies
$e_\ell=\{a',\hat{I}_{{\phi_{a'}}_-}(a)\}$ in the $\hat{I}_t$ process,
and the latter must be in the event $R_b$. This shows that
$\mathbb{P}(L_b)\leq\mathbb{P}(R_b)$. The opposite inequality
follows from reversing
roles of the two processes.
\end{pf*}

\subsection{\texorpdfstring{Estimate for the first depinking time (Lemma~\protect\ref{lemfirstdepinking})}
{Estimate for the first depinking time (Lemma 9.2)}}\label{secprooffirstdepinking}
\mbox{}
\begin{pf*}{Proof of Lemma~\protect\ref{lemfirstdepinking}} As in Lemma~\ref
{lemfirsttworounds}
we drop condition $(2)$ for a depinking time, and notice that this
change does \textit{not} change the value (or the distribution) of $D_1$.
The modification to the process also does not affect the end result of
Lemma~\ref{lemdiscreteMarkov}: that is, the discrete-time process starting
from $(\bz,R,P,W)$ with subsequent states
$(\hat{\bz}_i,\hat{R}_i,\hat{P}_i,\hat{W}_i)$ described in that lemma
is a time-homogeneous Markov chain, and $\hat{D}_1\equiv D_1/2T$ is a
stopping time for this process.

Moreover, we assume without loss of generality that $|R|\leq|W|$,
which implies that $|R^C_t|\leq|W^C_t|$ for $t<D_1$. This implies
$|W^C_t|\geq|V|/6$ unless there are more pink particles than red ones
at time $t<D_1$; this follows from the reasoning in the beginning of
the proof of Lemma~\ref{lemfirsttworounds} in Section~\ref
{secprooffirsttworounds}.
Recall that each pinkening step remores a red particle and creates two
pink ones. It follows that $|R^C_{2iT_-}|<2|R|/3$ implies $D_1\leq
2iT$, and
%
%
\setcounter{equation}{0}
\begin{eqnarray}\forall i\in\N_+\qquad  \mathbb{P}(\hat
{D}_1>i)&\leq
&\mathbb{P}(|\hat{R}_i|\geq2|R|/3, \hat{D}_1>i-1)
\nonumber
\\[-8pt]
\\[-8pt]
\nonumber
&\leq&\frac
{3\mathbb
{E}[\hat{R}_{i}{\mathbb{I}_{\{\hat{D}_1>i\}}}]}{2|R|}\\
\label
{eqintegrandoBURRO!}&=&
\frac{3\mathbb{E}[\mathbb{E}[\hat{R}_{i}\vert\hat{\sF
}_{i-1}]{\mathbb{I}_{\{\hat{D}_1>(i-1)\}}}]}{2|R|},
\end{eqnarray}
where $\hat{\sF}_{i-1}$ is the $\sigma$-field generated by $(\hat
{\bz
}_\ell,\hat{R}_\ell,\hat{P}_\ell,\hat{W}_\ell)$ for $\ell\leq i-1$.

We now estimate the integrand in (\ref{eqintegrandoBURRO!}).
Lemma~\ref{lemdiscreteMarkov} and its proof impliy that
\[
\mathbb{E}[|\hat{R}_{i}|\vert\hat{\sF}_{i-1}]
\]
is the expected number of red particles after a potential depinking, a
constant-color phase and a color-changing phase for a chameleon process
started from
\[
(\hat{\bz}_{i-1},\hat{R}_{i-1},\hat{P}_{i-1},\hat{W}_{i-1})\in\sC_k(V).
\]
By Lemma~\ref{lemfirsttworounds}, we can ensure that
\[
\mathbb{E}[|\hat{R}_{i}|\vert\hat{\sF}_{i-1}]\leq
(1-c)|\hat{R}_{i-1}|\qquad\mbox{if
}|\hat{W}_{i-1}|\geq|V|/6\quad\mbox{and}\quad|\hat{P}_{i-1}|<|\hat{R}_{i-1}|.
\]
As noted before, these conditions are always satisfied in the event $\{
\hat{D}_1>(i-1)\}$, because there are less pink than red particles. We deduce
\begin{eqnarray*}
\forall i\in\N_+\qquad  \frac{\mathbb
{E}[|\hat{R}_{i}|{\mathbb{I}_{\{\hat{D}_1>i\}}}]}{|R|}
&\leq&
\frac{\mathbb{E}[\mathbb{E}[|\hat{R}_{i}|\vert\hat
{\sF}_{i-1}]{\mathbb{I}_{\{\hat{D}_1>(i-1)\}}}
]}{|R|}\\
&\leq&
(1-c)\biggl\{\frac{\mathbb{E}[|\hat{R}_{i-1}|{\mathbb{I}_{\{
\hat {D}_1>i-1\} }}]}{|R|}\biggr\}\\
\mbox{(\ldots induction\ldots )}&\leq&(1-c)^i.
\end{eqnarray*}
This implies
\[
\mathbb{P}(D_1>2iT) =\mathbb{P}(\hat{D}_1>i)\leq\frac
{3(1-c)^i}{2},\qquad\mbox{$c=1/1000$
universal.}
\]
From this one can easily show that $\mathbb{E}[e^{D_1/KT}
]\leq e$ for some
universal $K$.
\end{pf*}

\subsection{\texorpdfstring{Proof of Lemma~\protect\ref{lemexpmomentDj}}{Proof of Lemma 6.2}}\label{secproofexpmomentDj}
\mbox{}
\begin{pf}Fix $\bx\in(V)_k$. We first prove that
%
%
\setcounter{equation}{0}
\begin{equation}\label{eqexpmomentnofill}\mathbb{E}\bigl[e^{
{D_j(\bx)}/{(K  T)}}\bigr]\leq e^j, \qquad\mbox{$K>0$ from Lemma~\ref
{lemfirstdepinking};}
\end{equation}
this is the bound we wish to obtain except that we are not conditioning
on $\Fill$.

We proceed as in the previous proof and consider the discrete-time process
\[
\{(\hat{\bz}_i,\hat{R}_i,\hat{P}_i,\hat{W}_i)\}_{i\geq0},
\]
introduced in Lemma~\ref{lemdiscreteMarkov}, henceforth called the hat
process. This time we take the initial state
\[
(\bz,R,P,W)\equiv\bigl(\bz,\{x\},\varnothing,V\setminus\bigl(\bO(\bz)\cup\{
x\}\bigr)\bigr)\vadjust{\goodbreak}
\]
corresponding to $\bx=(\bz,x)$ in the sense of
Proposition~\ref{propchameleonisrelated}. Also recall the definition
$\hat{D}_i\equiv D_{i}(\bx)/2T$ and note that
(\ref{eqexpmomentnofill}) is equivalent to
%
%
\begin{equation}\label{eqexpmomentnofill2}\mathbb{E}[e^{
{\hat{D}_j}/{K'}}]\leq e^j,\qquad \mbox{$K'=2K$.}
\end{equation}
This is valid for $j=1$ due to Lemma~\ref{lemfirstdepinking}. For $j>1$,
we recall the definition of the $\sigma$-fields $\hat{\sF}_i$,
recall that $\hat{D}_{j-1}$ is a stopping time for the hat process
(cf. Lemma~\ref{lemdiscreteMarkov}) and obtain
%
%
\begin{equation}\label{eqexpmomentnofill3}\mathbb{E}[e^{
{\hat{D}_j}/{K'}}]\leq\mathbb{E}\bigl[e^{{\hat
{D}_{j-1}}/{K'}} \mathbb{E}\bigl[e^{{(\hat{D}_{j}-\hat
{D}_{j-1})}/{K'}}\vert\hat{\sF}_{\hat{D}_{j-1}}\bigr]\bigr].
\end{equation}
We will apply the strong Markov property of the hat process (cf.
Lemma~\ref{lemdiscreteMarkov} again) to bound the conditional expectation
in the RHS. The conditional law of $\hat{D}_j-\hat{D}_{j-1}$ given
$\hat
{\sF}_{\hat{D}_{j-1}}$ is the law of the hat process
started from state
\[
(\hat{\bz}_{\hat{D}_{j-1}},\hat{R}_{\hat{D}_{j-1}},\hat{P}_{\hat
{D}_{j-1}},\hat{W}_{\hat{D}_{j-1}}).
\]
Notice that $\hat{P}_{\hat{D}_{j-1}} = P^C_{{D_{j-1}}_-}\neq
\varnothing$; in fact, since depinking occurs at time $D_{j-1}$, we
know that $|P^C_{{D_{j-1}}_-}|\geq
\min\{|R^C_{{D_{j-1}}_-}|,|W^C_{{D_{j-1}}_-}|\}$. However, at time
$D_{j-1}$ all pink particles disappear: the hat process evolves as
if started from a state with no pink particles, and
$\hat{D}_{j}-\hat{D}_{j-1}$ is the first depinking time for the hat
process with this modified initial state. We deduce from
Lemma~\ref{lemfirstdepinking} that
\[
\mathbb{E}\bigl[e^{{(\hat{D}_{j}-\hat{D}_{j-1})}/{K'}}\vert\hat
{\sF}_{\hat{D}_{j-1}}\bigr]\leq e\qquad \mbox{almost surely,}
\]
so that
%
%
\begin{equation}\label{eqexpmomentnofill4}\mathbb{E}[e^{
{\hat{D}_j}/{K'}}]\leq\mathbb{E}[e^{{\hat
{D}_{j-1}}/{K'}}] e\leq e^{j} \qquad\mbox{by induction.}
\end{equation}
This proves (\ref{eqexpmomentnofill3}) and (\ref{eqexpmomentnofill}).

To prove the lemma, notice that conditioning on $\Fill$ simply
biases the coin flips $d_i$ performed at depinking times; cf.
Remark~\ref{remproofdiscreteinkdecay}. This will not change the
distribution of $\hat{D}_1$ or the conditional distribution of
$\hat{D}_{j}-\hat{D}_{j-1}$ given the past of the process, so the
argument we presented above still applies.
\end{pf}
%
\section{Final remarks}\label{secconclusion}

Our paper leaves many questions open. Here we present a few problems
that seem especially interesting:

\begin{itemize}
\item Are there any other interacting particle systems whose mixing
parameters can be bounded solely in terms of the constituent
parts? Nonsymmetric exclusion is an obvious candidate. Another is
the zero-range process. Morris~\cite{MorrisZeroRange} used the
comparison principle and a
coupling argument on the complete graph to bound the spectral gap of
this process on a grid. Can one do something less indirect over an
arbitrary graph?
\item Can we find a mixing time upper bound of $\IP(|V|,G)$ (i.e., as
many particles as vertices), that is, similar to our main Theorem?
Inspection of the chameleon process shows that it gives the
conditional distribution of a particle given the \textit{whole past
trajectory} of the other particles. This means, in particular, that
it cannot deal with $k=|V|$ particles.
\item Recall the heuristic assumption in the \hyperref[sec1]{Introduction}:
${\fT}_{\EX(k,G)}(\eps)\leq C_1 \times\break{\fT}_{\RW(G)} (\eps/k)$ with $C_1>0$
universal. Is this actually
true? This would be stronger than our main theorem.
\item Combining the previous two items: is it true that
${\fT}_{\IP(|V|,G)}(\eps)=C_1 \times\break{\fT}_{\RW(G)} (\eps/|V|)$? Could it
even be
possible that ${\fT}_{\IP(|V|,G)}(\eps)\leq {\fT}_{\RW(|V|,G)}(\eps)$, that
is, the
interchange process mixes at least as fast as independent random
walkers? This would give Aldous's (now proven) conjecture on the
spectral gap as a corollary.
\end{itemize}

\begin{appendix}\label{app}

\section{Mixing bounds for $\EX(k,G)$ via canonical paths}\label
{secfaircomparison}

We use asymptotic notation below as shorthand; see, for example, \cite
{AlonSpencerMethod} for precise definitions. Let $G=(V,E,\{w_e\}_{e\in
E})$ be a weighted graph. It seems that the best \textit{general} bound that was previously available (implicitly) for the
mixing time of $\EX(k,G)$ comes from the combination of three
ingredients.

\textit{Mixing time from Log-Sobolev constant.} The state
space of $\EX(k,G)$ has cardinality ${{|V|}\choose{k}}= 2^{\Theta(|V|)}$
if $k=\Theta(|V|)$. By the results of~\cite{DiaconisSCLogSobolev},
if $\rho_{\EX(k,G)}$ is the log-Sobolev constant of $\EX(k,G)$,
then
\[
{\fT}_{\EX(k,G)}(1/4) = O\biggl(\frac{\ln|V|}{\rho_{\EX (k,G)}}\biggr)\qquad\mbox{for }k=\Theta(|V|).
\]

\textit{Log-Sobolev inequality for the Bernoulli--Laplace model.}
Consider the complete graph $K_{|V|}$ where each edge
has weight $1/|V|$. $\EX(k,K_{|V|})$ is the so-called \textit{Bernoulli--Laplace model} with
$k$ particles, whose log Sobolev constant is of the order $\Theta(\ln
({|V|}^2/k({|V|}-k)))$. Notice that this is $\Theta(1)$ for
$k=\Theta({|V|}).$

\textit{Comparison argument}. Now consider a general weighted graph
$G=(V,E,\break \{w_e\}_{e\in E})$. Assume that for
each pair $(x,y)\in V^2$ one has defined a path $\gamma_{x,y}$ in
$G$ connecting $x$ to $y$. For each such pair, let $I_{x,y}(e)=1$ if
$e$ is crossed by $\gamma_{x,y}$ and $0$ otherwise, and also let
$\ell_{x,y}$ denote the length of $\gamma_{x,y}$. Finally, define
\[
\phi(G)\equiv\max_{e\in E}\sum_{(x,y)\in V^2} \frac{\ell_{x,y}
I_{x,y}(e)}{|V| w_e}.
\]
It is shown in the proof of \cite[Theorem 3.1]{DiaconisSCComparison}
that this comparison constant for the Dirichlet forms of $\RW(G)$ can
be ``lifted'' with no loss to $\EX(k,G)$. The comparison principle for
the log Sobolev constant~\cite{DiaconisSCLogSobolev} implies $\rho
_{\EX
(k,G)}=\Omega(\phi^{-1}(G))$. We deduce
%
%
\renewcommand{\theequation}{A.\arabic{subsection}.\arabic{equation}}
\setcounter{equation}{0}
\begin{equation}\label{eqdiaconisscbound}{\fT}_{\EX(k,G)}(1/4) =
O(\phi(G) \ln|V|)\qquad\mbox{if }k=\Theta(|V|).
\end{equation}

It can be very hard to find good upper bounds on $\phi(G)$ in
general, but the general lower bound we will present implies that
%
%
\begin{equation}\label{eqphiisbad}\phi(G)\geq
\frac{2\overline{\mathrm{dist}^2}}{\overline{d}},
\end{equation}
where
$\overline{\mathrm{dist}^2}$ is the average over all $(x,y)\in V^2$ of
the square of the graph-theoretic distance between $x$ and $y$, and
$\overline{d}$ is the average (weighted) degree in $G$. Indeed, it
suffices to see that
\begin{eqnarray*}
\phi(G)&\geq&\sum_{e\in E}\frac{w_e}{\sum_{f\in
E}w_f}\biggl(\sum_{(x,y)\in V^2} \frac{I_{x,y}(e) \ell_{x,y}}{|V| w_e}
\biggr)\\
&=& \frac{1}{\sum_{f\in
E}w_f}\biggl[\sum_{(x,y)\in V^2}\biggl(\frac{\sum_{e\in
E}I_{x,y}(e)\ell_{x,y}}{|V|}\biggr)\biggr]\\
\mbox{[use $\sum_e I_{x,y}(e)=\ell_{x,y}$]}&=& \frac{1}{|V|^{-1}\sum_{f\in
E}w_f}\biggl[\sum_{(x,y)\in V^2}\frac{\ell_{x,y}^2}{|V|^2}\biggr]\\
\mbox{[use $\ell_{x,y}\geq\mathrm{dist}(x,y)$]}&\geq
&\frac{1}{|V|^{-1}\sum_{f\in E}w_f}\biggl[\sum_{(x,y)\in V^2}\frac
{\mathrm{dist}(x,y)^2}{|V|^2}\biggr]= \frac{2\overline{\mathrm{dist}^2}}{\overline{d}}.
\end{eqnarray*}

We note that this is a lower bound, which we do not know how to achieve
in the examples in Table~\ref{table1}.

\section{\texorpdfstring{The trajectory of $\widehat{\lowercase{\mathsf{ink}}_j}$ given $\mathsf{F}\lowercase{\mathsf{ill}}$}
{The trajectory of ink j given Fill}}\label{secproofdiscreteinkdecay}

We use the facts proven in Section~\ref{secmiscellany} to derive the
technical estimate in Proposition~\ref{propdiscreteinkdecay} in the
proof of
Lemma~\ref{lemkto2}; cf. Section~\ref{secproofkto2}.

\begin{pf*}{Proof of Proposition~\protect\ref{propdiscreteinkdecay}} We
take the notation
in Section~\ref{secmiscellany} for granted, but omit the superscript
$\bw$ in this proof. Our first goal will be to show that, conditionally
on $\Fill$,
$\{\hatink_{j}\}_{j\geq0}$ is still a Markov chain. Repeating the
steps of the proof of Lemma~\ref{lemsummaryhatink}, we note that
\[
\mathbb{P}(\Fill\vert(\hatink_{i})_{i\leq j}) =\frac{\mathbb
{E}
[\hatink_{\infty}\vert(\hatink_{i})_{i\leq j}
]}{|V|-k+1}=\frac{\hatink_{j}}{|V|-k+1} = \mathbb{P}(\Fill)
\hatink_j.
\]
We deduce
from Bayes's rule and the Markovian property that
\begin{eqnarray*}
\mathbb{P}\Biggl(\bigcap_{i=1}^{j}\{\hatink_i=a_i\}\vert
\Fill\Biggr)&=&
\mathbb{P}\Biggl(\bigcap_{i=1}^{j}\{\hatink_i=a_i\}\Biggr)  a_j\\ \mbox
{(Markov property
for $\hatink_j$)}&=& p(1,a_1)p(a_1,a_2)\cdots p(a_{j-1},a_j)  a_j\\
&=& q(1,a_1)\cdots q(a_{j-1},a_j),
\end{eqnarray*}
where
\[
q(a,b) = \frac{bp(a,b)}{a}\qquad\mbox{if $a\neq0$.}
\]
Notice that, since $\hatink_j$ does not visit $0$ in the event
$\Fill$, we do not need to define $q(a,b)$ for $a=0$. We have shown:

\begin{proposition}\label{propconditionalMarkov}Conditionally on
$\Fill
$, the trajectory of $\{\hatink_j\}_{j\geq
0}$ is that of a Markov chain in $\{1,\ldots,|V|-k+1\}$, with
transition rates $q(a,b)$ and started from
$\hatink_{0}=1$.
\end{proposition}

For the remainder of the proof, we will use this proposition to bound
$1-\hatink_\ell/(|V|-k+1)$. Actually, another quantity is
easier to bound. Set $I_\ell= \hatink_\ell/(|V|-k+1)$ and
\[
Z_\ell\equiv\frac{\sqrt{\min\{1-I_\ell,I_\ell\}}}{I_\ell}.
\]
Notice that conditionally on $\Fill$, $I_\ell>0$ always, hence
$Z_\ell$ is a.s. well defined for all $\ell$. Moreover, one can
check that $1-I_\ell\leq Z_\ell$ always. Therefore the lemma will
follow from the estimate
\[
\mathbb{E}^{\mathsf{Fill}}[Z_\ell]\leq(71/72)^\ell\sqrt{|V|-k+1},
\]
where $\mathbb{E}^{\mathsf{Fill}}[\cdot]$ corresponds to an
expectation with respect to
the conditional distribution given $\Fill$. Since $Z_0 =
\sqrt{|V|-k+1}$, the above estimate follows directly from the
following \textit{claim}.
\begin{claim}
\[
\forall\ell\in\N\qquad  \mathbb{E}^{\mathsf{Fill}}[Z_\ell]\leq
(71/72)
\mathbb{E}^{\mathsf{Fill}}[Z_{\ell-1}].
\]
\end{claim}

Therefore, proving this claim will
finish the proof.

To prove the claim, we first note that for all $i$, $Z_i$ is a
function of $\hatink_{i}$, and $Z_{\ell-1}=0\Rightarrow Z_\ell=0$.
We deduce
%
%
\renewcommand{\theequation}{B.\arabic{subsection}.\arabic{equation}}
\setcounter{equation}{0}
\begin{equation}\label{eqcovetedformula}\mathbb{E}^{\mathsf
{Fill}}[Z_\ell] =
\mathbb{E}^{\mathsf{Fill}}\biggl[\mathbb{E}^{\mathsf{Fill}}\biggl[\frac
{Z_{\ell}}{Z_{\ell-1}}\Bigm| \hatink_{\ell-1}\biggr] Z_{\ell-1}{\mathbb
{I}_{\{Z_{\ell-1}\neq 0\}}}\biggr].
\end{equation}
We now bound the conditional expectation in the RHS. We
may assume that $\hatink_{\ell-1}=r$ with $0<r<|V|-k+1$ (otherwise
$Z_{\ell-1}=0$). Thus we wish to bound
\[
\mathbb{E}^{\mathsf{Fill}}\biggl[\frac{Z_{\ell}}{Z_{\ell-1}}\Bigm|\hatink
_{\ell-1}=r\biggr],\qquad
1\leq
r\leq|V|-k.
\]
If we note that
\[
\frac{Z_\ell}{Z_{\ell-1}} =
\frac{\sqrt{\min\{1-I_\ell,I_\ell\}}}{\sqrt{\min\{1-I_{\ell
-1},I_{\ell
-1}\}}}\times
\frac{I_{\ell-1}}{I_\ell} = \frac{\sqrt{\min\{1-I_\ell,I_\ell\}
}}{\sqrt
{\min\{1-I_{\ell-1},I_{\ell-1}\}}}\times
\frac{\hatink_{\ell-1}}{\hatink_\ell}
\]
and define $f(a) = \sqrt{\min\{a,|V|-k+1-a\}}$, we see that
\begin{eqnarray*}
\mathbb{E}^{\mathsf{Fill}}\biggl[\frac{Z_{\ell}}{Z_{\ell
-1}}\Big\vert\hatink _{\ell-1}=r\biggr]
&=& \mathbb{E}^{\mathsf{Fill}}\biggl[\frac{f(\hatink_{\ell})}{f(\hatink
_{\ell-1})}\times \frac{\hatink_{\ell-1}}{\hatink_\ell}\Big\vert
\hatink_{\ell-1}=r\biggr]\\
\mbox{(use Proposition~\ref{propconditionalMarkov})}&=&
\sum_{s}q(r,s) \frac{f(s)}{f(r)}\times\frac{r}{s}\\
\mbox{[use formula for $q(\cdot,\cdot\cdot)$]} &=& \sum_{s}p(r,s)
\frac{f(s)}{f(r)}
\end{eqnarray*}
where $p(\cdot,\cdot\cdot)$ are the
transition rates of the unconditional $\{\hatink_j\}_{j\geq0}$
process. Using the formula for these, we obtain
%
%
\begin{equation}\label{eqcovetedconcave}\mathbb{E}^{\mathsf
{Fill}}\biggl[\frac{Z_{\ell }}{Z_{\ell -1}}\Bigm|\hatink_{\ell-1}=r\biggr] =
\frac{1}{2}\biggl(\frac{f(r+\Delta(r)) +
f(r-\Delta(r))}{f(r)}\biggr).
\end{equation}

Recall the formula for $\Delta(r)$ (cf. Proposition~\ref{propinkincrement}),
\[
\Delta(r) \equiv
\biggl\lceil\frac{\min\{r,|V|-k+1-r\}}{3}\biggr\rceil.
\]

We now split the analysis of the RHS of this in two cases.

\textit{Case} 1: $1\leq r\leq(|V|-k+1)/2$. In this case
$f(r)=\sqrt{r}$ and $\Delta(r) = \lceil r/3\rceil\geq r/3$. We use
the upper bound $f(r\pm\Delta(r))\leq\sqrt{r\pm\Delta(r)}$ to
obtain
%
%
\begin{equation}\label{eqcovetedconcave1i}\mathbb{E}^{\mathsf
{Fill}}\biggl[\frac{Z_{\ell }}{Z_{\ell-1}}\Bigm|\hatink_{\ell-1}=r\biggr] =
\frac{1}{2}\Biggl(\sqrt{1-\frac{\Delta(r)}{r}} +
\sqrt{1+\frac{\Delta(r)}{r}}\Biggr).
\end{equation}
Recall the bound
``$\sqrt{1-x}+\sqrt{1+x}\leq2(1-x^2/8)$,'' valid for all $0\leq
x\leq1$; this can be checked by squaring both sides of the
inequality. In our case, we apply this with $x=\Delta(r)/r\geq1/3$
and deduce
%
%
\begin{equation}\label{eqcovetedconcave1}\mathbb{E}^{\mathsf
{Fill}}\biggl[\frac{Z_{\ell }}{Z_{\ell-1}}\Bigm|\hatink_{\ell-1}=r\biggr] =
1 - \frac{1}{8} \biggl(\frac{\Delta(r)}{r}\biggr)^2\leq
\frac{71}{72}.
\end{equation}

\textit{Case} 2: $(|V|-k+1)/2<r\leq|V|-k$. In this case (\ref
{eqcovetedconcave1i}) holds with $r'=|V|-k+1-r$ replacing $r$. Similar
calculations imply that the conditional expectation is also $\leq
71/72$ in this case.

Thus we see that in both cases
\[
\mathbb{E}^{\mathsf{Fill}}\biggl[\frac{Z_{\ell}}{Z_{\ell-1}}\Bigm|\hatink
_{\ell-1}=r\biggr]\leq
\frac{71}{72}.
\]

Plugging this into (\ref{eqcovetedformula}) gives
\[
\mathbb{E}^{\mathsf{Fill}}[Z_\ell]\leq\frac{71}{72}  \mathbb
{E}^{\mathsf{Fill}}\bigl[Z_{\ell-1}{\mathbb{I}_{\{ Z_{\ell -1}\neq0\}
}}\bigr] = \frac{71}{72}  \mathbb{E}^{\mathsf{Fill}}[Z_{\ell-1}],
\]
which completes the proof.
\end{pf*}
\end{appendix}
\section*{Acknowledgments}
We thank Ton Dieker and Prasad Tetali for
useful discussions on the exposition. We also thank an anonymous
referee for a long list of typos in the previous version, as well as
for numerous suggestions.
%



\printaddresses


\begin{thebibliography}{24}

\bibitem{AldousFillRWBook}
\begin{bmisc}[auto:STB|2012/02/29|12:31:17]
\bauthor{\bsnm{Aldous},~\bfnm{David}\binits{D.}} \AND
  \bauthor{\bsnm{Fill},~\bfnm{James~Allen}\binits{J.~A.}}
\bhowpublished{Reversible Markov {Chains} and random {walks} on graphs. Book
  draft. Available at \texttt{%
  \href{http://www.stat.berkeley.edu/\textasciitilde aldous/RWG/book.html}{http://www.stat.berkeley.edu/\textasciitilde aldous/}
  \href{http://www.stat.berkeley.edu/\textasciitilde aldous/RWG/book.html}{RWG/book.html}}}.
\bptok{imsref}%
\end{bmisc}
\endbibitem

\bibitem{AlonSpencerMethod}
\begin{bbook}[mr]
\bauthor{\bsnm{Alon},~\bfnm{Noga}\binits{N.}} \AND
  \bauthor{\bsnm{Spencer},~\bfnm{Joel~H.}\binits{J.~H.}}
(\byear{2000}).
\btitle{The Probabilistic Method},
\bedition{2nd} ed.
\bpublisher{Wiley}, \baddress{New York}.
\bid{doi={10.1002/0471722154}, mr={1885388}}
\bptok{imsref}%
\end{bbook}
\endbibitem

\bibitem{AndjelCorrelationSymmetricExclusion}
\begin{barticle}[mr]
\bauthor{\bsnm{Andjel},~\bfnm{Enrique~D.}\binits{E.~D.}}
(\byear{1988}).
\btitle{A correlation inequality for the symmetric exclusion process}.
\bjournal{Ann. Probab.}
\bvolume{16}
\bpages{717--721}.
\bid{issn={0091-1798}, mr={0929073}}
\bptok{imsref}%
\end{barticle}
\endbibitem

\bibitem{BenjaminiMosselRWPercolation}
\begin{barticle}[mr]
\bauthor{\bsnm{Benjamini},~\bfnm{Itai}\binits{I.}} \AND
  \bauthor{\bsnm{Mossel},~\bfnm{Elchanan}\binits{E.}}
(\byear{2003}).
\btitle{On the mixing time of a simple random walk on the super critical
  percolation cluster}.
\bjournal{Probab. Theory Related Fields}
\bvolume{125}
\bpages{408--420}.
\bid{doi={10.1007/s00440-002-0246-y}, issn={0178-8051}, mr={1967022}}
\bptok{imsref}%
\end{barticle}
\endbibitem

\bibitem{CaputoFaggionatoRWPointProcess}
\begin{barticle}[mr]
\bauthor{\bsnm{Caputo},~\bfnm{Pietro}\binits{P.}} \AND
  \bauthor{\bsnm{Faggionato},~\bfnm{Alessandra}\binits{A.}}
(\byear{2007}).
\btitle{Isoperimetric inequalities and mixing time for a random walk on a
  random point process}.
\bjournal{Ann. Appl. Probab.}
\bvolume{17}
\bpages{1707--1744}.
\bid{doi={10.1214/07-AAP442}, issn={1050-5164}, mr={2358639}}
\bptok{imsref}%
\end{barticle}
\endbibitem

\bibitem{CaputoEtAlInterchange}
\begin{barticle}[mr]
\bauthor{\bsnm{Caputo},~\bfnm{Pietro}\binits{P.}},
  \bauthor{\bsnm{Liggett},~\bfnm{Thomas~M.}\binits{T.~M.}} \AND
  \bauthor{\bsnm{Richthammer},~\bfnm{Thomas}\binits{T.}}
(\byear{2010}).
\btitle{Proof of {A}ldous' spectral gap conjecture}.
\bjournal{J. Amer. Math. Soc.}
\bvolume{23}
\bpages{831--851}.
\bid{doi={10.1090/S0894-0347-10-00659-4}, issn={0894-0347}, mr={2629990}}
\bptok{imsref}%
\end{barticle}
\endbibitem

\bibitem{CooperEtAlMultipleRWIPS}
\begin{bincollection}[mr]
\bauthor{\bsnm{Cooper},~\bfnm{Colin}\binits{C.}},
  \bauthor{\bsnm{Frieze},~\bfnm{Alan}\binits{A.}} \AND
  \bauthor{\bsnm{Radzik},~\bfnm{Tomasz}\binits{T.}}
(\byear{2009}).
\btitle{Multiple random walks and interacting particle systems}.
In \bbooktitle{Automata, Languages and Programming. {P}art {II}}.
\bseries{Lecture Notes in Computer Science}
\bvolume{5556}
\bpages{399--410}.
\bpublisher{Springer}, \baddress{Berlin}.
\bid{mr={2544812}}
\bptok{imsref}%
\end{bincollection}
\endbibitem

\bibitem{DiaconisSCComparison}
\begin{barticle}[mr]
\bauthor{\bsnm{Diaconis},~\bfnm{Persi}\binits{P.}} \AND
  \bauthor{\bsnm{Saloff-Coste},~\bfnm{Laurent}\binits{L.}}
(\byear{1993}).
\btitle{Comparison theorems for reversible {M}arkov chains}.
\bjournal{Ann. Appl. Probab.}
\bvolume{3}
\bpages{696--730}.
\bid{issn={1050-5164}, mr={1233621}}
\bptok{imsref}%
\end{barticle}
\endbibitem

\bibitem{DiaconisSCLogSobolev}
\begin{barticle}[mr]
\bauthor{\bsnm{Diaconis},~\bfnm{P.}\binits{P.}} \AND
  \bauthor{\bsnm{Saloff-Coste},~\bfnm{L.}\binits{L.}}
(\byear{1996}).
\btitle{Logarithmic {S}obolev inequalities for finite {M}arkov chains}.
\bjournal{Ann. Appl. Probab.}
\bvolume{6}
\bpages{695--750}.
\bid{doi={10.1214/aoap/1034968224}, issn={1050-5164}, mr={1410112}}
\bptok{imsref}%
\end{barticle}
\endbibitem

\bibitem{DiekerInterchange}
\begin{barticle}[mr]
\bauthor{\bsnm{Dieker},~\bfnm{A.~B.}\binits{A.~B.}}
(\byear{2010}).
\btitle{Interlacings for random walks on weighted graphs and the interchange
  process}.
\bjournal{SIAM J. Discrete Math.}
\bvolume{24}
\bpages{191--206}.
\bid{doi={10.1137/090775361}, issn={0895-4801}, mr={2600660}}
\bptok{imsref}%
\end{barticle}
\endbibitem

\bibitem{FountoulakisReedMixingGnp}
\begin{barticle}[mr]
\bauthor{\bsnm{Fountoulakis},~\bfnm{N.}\binits{N.}} \AND
  \bauthor{\bsnm{Reed},~\bfnm{B.~A.}\binits{B.~A.}}
(\byear{2008}).
\btitle{The evolution of the mixing rate of a simple random walk on the giant
  component of a random graph}.
\bjournal{Random Structures Algorithms}
\bvolume{33}
\bpages{68--86}.
\bid{doi={10.1002/rsa.20210}, issn={1042-9832}, mr={2428978}}
\bptok{imsref}%
\end{barticle}
\endbibitem

\bibitem{LeeYauLogSobolev}
\begin{barticle}[mr]
\bauthor{\bsnm{Lee},~\bfnm{Tzong-Yow}\binits{T.-Y.}} \AND
  \bauthor{\bsnm{Yau},~\bfnm{Horng-Tzer}\binits{H.-T.}}
(\byear{1998}).
\btitle{Logarithmic {S}obolev inequality for some models of random walks}.
\bjournal{Ann. Probab.}
\bvolume{26}
\bpages{1855--1873}.
\bid{doi={10.1214/aop/1022855885}, issn={0091-1798}, mr={1675008}}
\bptok{imsref}%
\end{barticle}
\endbibitem

\bibitem{LevinPeresWilmerMCBook}
\begin{bbook}[mr]
\bauthor{\bsnm{Levin},~\bfnm{David~A.}\binits{D.~A.}},
  \bauthor{\bsnm{Peres},~\bfnm{Yuval}\binits{Y.}} \AND
  \bauthor{\bsnm{Wilmer},~\bfnm{Elizabeth~L.}\binits{E.~L.}}
(\byear{2009}).
\btitle{Markov Chains and Mixing Times}.
\bpublisher{Amer. Math. Soc.}, \baddress{Providence, RI}.
\bid{mr={2466937}}
\bptnote{check year}%
\bptok{imsref}%
\end{bbook}
\endbibitem

\bibitem{LiggettInvariantMeasures2}
\begin{barticle}[mr]
\bauthor{\bsnm{Liggett},~\bfnm{Thomas~M.}\binits{T.~M.}}
(\byear{1974}).
\btitle{A characterization of the invariant measures for an infinite particle
  system with interactions. {II}}.
\bjournal{Trans. Amer. Math. Soc.}
\bvolume{198}
\bpages{201--213}.
\bid{issn={0002-9947}, mr={0375531}}
\bptok{imsref}%
\end{barticle}
\endbibitem

\bibitem{LiggettIPSBook}
\begin{bbook}[mr]
\bauthor{\bsnm{Liggett},~\bfnm{Thomas~M.}\binits{T.~M.}}
(\byear{1985}).
\btitle{Interacting Particle Systems}.
\bseries{Grundlehren der Mathematischen Wissenschaften [Fundamental Principles
  of Mathematical Sciences]}
\bvolume{276}.
\bpublisher{Springer}, \baddress{New York}.
\bid{mr={0776231}}
\bptok{imsref}%
\end{bbook}
\endbibitem

\bibitem{LiggettSISBook}
\begin{bbook}[mr]
\bauthor{\bsnm{Liggett},~\bfnm{Thomas~M.}\binits{T.~M.}}
(\byear{1999}).
\btitle{Stochastic Interacting Systems: Contact, Voter and Exclusion
  Processes}.
\bseries{Grundlehren der Mathematischen Wissenschaften [Fundamental Principles
  of Mathematical Sciences]}
\bvolume{324}.
\bpublisher{Springer}, \baddress{Berlin}.
\bid{mr={1717346}}
\bptok{imsref}%
\end{bbook}
\endbibitem

\bibitem{MontenegroTetaliMixingBook}
\begin{barticle}[mr]
\bauthor{\bsnm{Montenegro},~\bfnm{Ravi}\binits{R.}} \AND
  \bauthor{\bsnm{Tetali},~\bfnm{Prasad}\binits{P.}}
(\byear{2006}).
\btitle{Mathematical aspects of mixing times in {M}arkov chains}.
\bjournal{Found. Trends Theor. Comput. Sci.}
\bvolume{1}
\bpages{x+121}.
\bid{issn={1551-305X}, mr={2341319}}
\bptok{imsref}%
\end{barticle}
\endbibitem

\bibitem{MorrisZeroRange}
\begin{barticle}[mr]
\bauthor{\bsnm{Morris},~\bfnm{Ben}\binits{B.}}
(\byear{2006}).
\btitle{Spectral gap for the zero range process with constant rate}.
\bjournal{Ann. Probab.}
\bvolume{34}
\bpages{1645--1664}.
\bid{doi={10.1214/009117906000000304}, issn={0091-1798}, mr={2271475}}
\bptok{imsref}%
\end{barticle}
\endbibitem

\bibitem{MorrisSimpleExclusion}
\begin{barticle}[mr]
\bauthor{\bsnm{Morris},~\bfnm{Ben}\binits{B.}}
(\byear{2006}).
\btitle{The mixing time for simple exclusion}.
\bjournal{Ann. Appl. Probab.}
\bvolume{16}
\bpages{615--635}.
\bid{doi={10.1214/105051605000000728}, issn={1050-5164}, mr={2244427}}
\bptok{imsref}%
\end{barticle}
\endbibitem

\bibitem{MorrisImproved}
\begin{barticle}[mr]
\bauthor{\bsnm{Morris},~\bfnm{Ben}\binits{B.}}
(\byear{2009}).
\btitle{Improved mixing time bounds for the {T}horp shuffle and {$L$}-rever\-sal
  chain}.
\bjournal{Ann. Probab.}
\bvolume{37}
\bpages{453--477}.
\bid{doi={10.1214/08-AOP409}, issn={0091-1798}, mr={2510013}}
\bptok{imsref}%
\end{barticle}
\endbibitem

\bibitem{MorrisPeresEvolvingSets}
\begin{barticle}[mr]
\bauthor{\bsnm{Morris},~\bfnm{B.}\binits{B.}} \AND
  \bauthor{\bsnm{Peres},~\bfnm{Yuval}\binits{Y.}}
(\byear{2005}).
\btitle{Evolving sets, mixing and heat kernel bounds}.
\bjournal{Probab. Theory Related Fields}
\bvolume{133}
\bpages{245--266}.
\bid{doi={10.1007/s00440-005-0434-7}, issn={0178-8051}, mr={2198701}}
\bptok{imsref}%
\end{barticle}
\endbibitem

\bibitem{PeteConnectivityMixingPercolation}
\begin{barticle}[mr]
\bauthor{\bsnm{Pete},~\bfnm{G{\'a}bor}\binits{G.}}
(\byear{2008}).
\btitle{A note on percolation on {$\Bbb Z\sp d$}: Isoperimetric profile via
  exponential cluster repulsion}.
\bjournal{Electron. Commun. Probab.}
\bvolume{13}
\bpages{377--392}.
\bid{doi={10.1214/ECP.v13-1390}, issn={1083-589X}, mr={2415145}}
\bptok{imsref}%
\end{barticle}
\endbibitem

\bibitem{YauLogSobolevExclusion}
\begin{barticle}[mr]
\bauthor{\bsnm{Yau},~\bfnm{Horng-Tzer}\binits{H.-T.}}
(\byear{1997}).
\btitle{Logarithmic {S}obolev inequality for generalized simple exclusion
  processes}.
\bjournal{Probab. Theory Related Fields}
\bvolume{109}
\bpages{507--538}.
\bid{doi={10.1007/s004400050140}, issn={0178-8051}, mr={1483598}}
\bptok{imsref}%
\end{barticle}
\endbibitem

\end{thebibliography}
\end{document}